\documentclass[a4paper, 12pt]{article}
\usepackage{xy}
\usepackage{epsfig}
\usepackage[T1]{fontenc} 
\usepackage[latin1]{inputenc}
\usepackage{amsfonts,amssymb,amsmath}
\usepackage{hyperref} 
\usepackage{graphics} 

\newtheorem{theo}{\sc Theorem}[section]

\newtheorem{defi}[theo]{\sc Definition}

\newtheorem{exam}[theo]{\sc Example}

\def\qed{ \ \hfil$\square$}

\def\Ker{{\hskip0.3mm\rm  Ker\hskip0.5mm}}
\def\Im{{\hskip0.3mm\rm  Im\hskip0.5mm}}

\makeindex

\newcommand{\Mt}{\widetilde{\mathcal M}}

\newcommand{\di}[1]{\displaystyle{#1}}

\newcommand{\Kpi}{\left(\frac{-1}{4\pi}\right)}

\newcommand{\pa}{\partial}

\newcommand{\tpa}{\tilde{\partial}}
\newcommand{\tr}{{\rm tr}}

\newcommand{\dd}{{\rm d}}

\newcommand{\GL}{{\rm GL}}

\newcommand{\Hom}{{\rm Hom}}

\newcommand{\M}{{{\rm M}}}

\newcommand{\ord}{{\rm ord}}

\newcommand{\tors}{{\rm tors}}

\newcommand{\lcm}{{\rm LCM~\!\!}}
\newcommand{\N}{{\mathbb N}}
\newcommand{\Q}{{\mathbb Q}}
\newcommand{\R}{{\mathbb R}}
\newcommand{\Z}{{\mathbb Z}}

\newcommand{\Qb}{\overline{\mathbb Q}}

\newcommand{\Ar}{{\mathcal A}}

\newcommand{\Br}{{\mathcal B}}

\newcommand{\Dr}{{\mathcal D}}

\newcommand{\Mc}{{\mathcal M}}
\newcommand{\Mr}{{\mathcal M}}
\newcommand{\Lr}{{\mathcal L}}

\newcommand{\Sr}{{\mathcal S}}

\newcommand{\al}{\alpha}
\newcommand{\De}{\Delta}
\newcommand{\de}{\delta}

\newcommand{\Ga}{\Gamma}
\newcommand{\ga}{\gamma}

\newcommand{\ka}{\kappa}
\newcommand{\La}{\Lambda}

\newcommand{\fragam}[2]{\frac{\Gamma_m(#1)}{\Gamma_m(#2)}}

\font\teneusm=eusm10 \font\seveneusm=eusm7 
\font\fiveeusm=eusm5 
\newfam\eusmfam 
\textfont\eusmfam=\teneusm 
\scriptfont\eusmfam=\seveneusm 
\scriptscriptfont\eusmfam=\fiveeusm

\def\mat #1,#2,#3,#4,{\left({#1\atop #3}{#2\atop #4}\right)}
\def\bra#1,{{\left\lbrace {#1}\right\rbrace}}

\def\si{\sigma}

\def \M{{\rm M}}

\def \Sp{{\rm Sp}}

\def\l1{\langle}

\newcommand{\B}{\left(\begin{array}{cc}}
\newcommand{\E}{\end{array}\right)}

\def \fns{{${}^{*)}$}}

\newcommand{\comm}[1]
{\fns\marginpar{$\boxed
{\hskip-6pt
{\small {\sf 
\begin{tabular} {l}
 #1
\end{tabular}
}
}
}
$
}
}
\def \?  {\comm{check ?}}

\usepackage{euscript}
\let\scr=\EuScript
\def\scrA{{\scr A}}   
\def\scrC{{\scr C}}

   \def\scrL{{\scr L}}
\def\scrM{{\scr M}}   
   \def\scrP{{\scr P}}
   
   \def\scrT{{\scr T}}

\def\Mr{{\scr M}}

\let\mathcal=\scr           
\def\ang#1,{{\left\langle {#1}\right\rangle}} 
\newcommand{\ds}{\displaystyle}

\def \rem  {\medskip\noindent {\sc Remark}}

\def\qed{ \ \hfil$\square$}

\def\Ker{{\hskip0.3mm\rm  Ker\hskip0.5mm}}
\def\Im{{\hskip0.3mm\rm  Im\hskip0.5mm}}

\font\teneusm=eusm10 \font\seveneusm=eusm7 
\font\fiveeusm=eusm5 
\newfam\eusmfam 
\textfont\eusmfam=\teneusm 
\scriptfont\eusmfam=\seveneusm 
\scriptscriptfont\eusmfam=\fiveeusm

\font\tengothic=eufm10
\font\sevengothic=eufm7
\font\fivegothic=eufm5
\newfam\Gothic
\textfont\Gothic=\tengothic
\scriptfont\Gothic=\sevengothic
\scriptscriptfont\Gothic=\fivegothic

\def\mat #1,#2,#3,#4,{\left({#1\atop #3}{#2\atop #4}\right)}
\def\bra#1,{{\left\lbrace {#1}\right\rbrace}}

\def\si{\sigma}

\def \M{{\rm M}}

\def \Sp{{\rm Sp}}

\def\l1{\langle}

\def\Numero{${\rm N\sp\circ}$}

\let\scr=\EuScript
\def\scrA{{\scr A}}   
\def\scrC{{\scr C}}

   \def\scrL{{\scr L}}
\def\scrM{{\scr M}}   
   \def\scrP{{\scr P}}
   
   \def\scrT{{\scr T}}

\def\Mr{{\scr M}}

\let\mathcal=\scr           
\let\db=\mathbb
   
\def\dbC{{\db C}}   
   
   \def\dbH{{\db H}}

   \def\dbR{{\db R}}

   \def\dbZ{{\db Z}}
\let\bold=\boldsymbol

\def\bpsi{{\bold \psi}}

\def\im{\Im}

\def\vin{{ {\tiny \mid }  
\kern-7.29pt 
\bigcup }}

\def\ang#1,{{\left\langle {#1}\right\rangle}} 

\def \cds{{\cdot{\dots}\cdot}}

\normalsize

\newcommand{\HH}{\mathbb H}

\newcommand{\ZZ}{\mathbb Z}

\def\cb{\C}

\def\zb{\Z}

\def\mc{\Mr}

\newcounter{ncours}{\setcounter{ncours} {1}}

\newcommand{\kp}{\kappa}
\newcommand{\vkp}{\varkappa}

\def\Lb{\Lambda}

           \def\cb{{\Bbb C}}

\def\mc{{\cal M}}

           \def\zb{{\Bbb Z}}

\def\Sp {\mathop{\rm Sp}\nolimits}

\title
{
 $p$-adic Banach modules of arithmetical modular forms\\
and triple products of Coleman's families
}
\author
{A. A.  Panchishkin \\ \noindent
{\it To dear Jean-Pierre Serre for his eightieth birthday with admiration}
\\
  Institut Fourier, 
  Université Grenoble-1\\
}
\date
{
Based on  a talk
  for the French-German Seminar  
in Lille 
on~November~9,~2005}

\def\ZZ{\mathbb Z}

\begin{document}
\maketitle

\begin{abstract}

For a prime number $p\ge 5$,
we consider three classical cusp eigenforms  
\begin{align*}
 f_j(z)=\sum_{n=1}^\infty a_{n,j}e(nz)\in \Sr_{k_j}(N_j, \psi_j),\  (j=1, 2,3)
\end{align*}
of weights  $k_1, k_2, k_3$, of 
 conductors $N_1, N_2, N_3$,
 and of nebentypus characters  $\psi_j \bmod N_j$.

According to H.Hida \cite{Hi86} and R.Coleman \cite{CoPB},
 one can include each $f_j$ $(j=1, 2, 3)$
 (under suitable  assumptions on $p$ and on $f_j$)  into a
 {$p$-adic analytic family}
{
$$
k_j{}\mapsto \{f_{j,k_j{}}= \sum_{n=1}^\infty a_{n}(f_{j, k_j{}})q^n\}
$$ }
{}
 of cusp eigenforms $f_{j,k_j{}}$ of weights $k_j{}$ 
 in such a way that $f_{j,k_j}=f_j$,
  and that all their Fourier coefficients $a_n(f_{j, k_j{}})$
  are given by 
 certain $p$-adic analytic functions
$k_j{}\mapsto a_{n, j}(k_j{})$.


The purpose of this paper is to describe a {\sf{} four variable}
 $p$-adic
$L$ function attached to Garrett's triple product of three Coleman's families
{}
$$
k_j{}\mapsto \left \{f_{j,k_j{}}= \sum_{n=1}^\infty a_{n,j}(k{})
q^n\right \} 
$$  
{}

 of cusp eigenforms of three fixed slopes $\sigma_j=v_p(\alpha_{p, j}^{(1)}(k_j{}))\ge 0$
 where $\alpha_{p,j}^{(1)} = \al_{p,j}^{(1)}(k_j{})$
 is an eigenvalue (which depends on $k_j{}$) of  Atkin's   operator $U=U_p$
 acting on Fourier expansions by
 $U(\sum_{n\ge 0}^\infty a_{n}q^n ) = \sum_{n \ge 0}^\infty a_{np} q^n$.


Let us consider the product of three eigenvalues:
$$
\lambda= \lambda(k_1{}, k_2{}, k_3{})
= \alpha_{p, 1}^{(1)}(k{}_1)\alpha_{p, 2}^{(1)}(k{}_2)\alpha_{p, 3}^{(1)}(k_3{})
$$
and assume that the slope of this product  
$$
\sigma= v_p(\lambda(k_1{}, k_2{}, k_3{}))=
\sigma(k_1{}, k_2{}, k_3{})=
\sigma_1+\sigma_2+\sigma_3
$$
 is {\sf constant and positive} for all triplets $(k_1{}, k_2{}, k_3{})$ 
 in an appropriate $p$-adic neighbourhood
 of the fixed triplet of weights $(k_1, k_2, k_3)$.
The each value  $\sigma_j$ is fixed.



We consider the $p$-adic weight space $X$ containing
all $(k{}_j, \psi_j)$.
Our $p$-adic $L$-functions are  {\sf Mellin transforms} of certain
 measures with values in $\Ar$, 
where $\Ar=\Ar({\cal  B})$ denotes  an affinoid algebra associated
 with an affinoid space ${\cal  B}$
 as in \cite{CoPB}, where 
${\cal  B}={\cal  B}_1\times{\cal  B}_2\times{\cal  B}_3,
$
 is an affinoid neighbourhood
 around 
$(k_1, k_2, k_3)\in X^3$ (with
 a given  integers $k_j$ and fixed Dirichlet characters 
$\psi_j \bmod N$).

We construct such a measure  from {\sf {}
 higher twists of classical
 Siegel-Eisenstein series}, 
 which produce distributions with values in
 certain Banach $\Ar$-modules $\Mr = \Mr(N;\Ar)$
 of triple modular forms with coefficients in the algebra 
$\Ar$.
\end{abstract}

  \tableofcontents

\section{Introduction}

\subsubsection*{Why study $L$-values attached to modular forms?}
{{} A popular proceedure in Number Theory is the following: }

\hskip-1cm
{}
\begin{tabular}{cccc} 
\begin{tabular}{l} 
Construct a generating \\ function
 $f=\sum_{n=0}^\infty a_n q^n$ \\
$\in {\mathbb C}[[ q]]$
of an arithmetical \\
function
$n\mapsto a_n$, 
\\
for example $a_n=p(n)$
\end{tabular}
& {}
$\rightsquigarrow$
\begin{tabular}{l}
Compute $f$ 
via \\
 modular forms, 
\\
for example
\\
$\ds\sum_{n=0}^\infty p(n) q^n$
\\
$=(\Delta/q)^{-1/24}$
\end{tabular}
&

$\rightsquigarrow$ 
\begin{tabular}{l}
A number
\\
(solution)
\end{tabular}
\\ 

\begin{tabular}{l}
Example 1 \rm \cite{Chand70}: 
 \\
(Hardy-Ramanujan)
\end{tabular} 
&

$\uparrow$
&
$\uparrow$
\\

\tiny
\begin{tabular}{l}
$\ds
 p(n)=\frac{e^{\pi \sqrt{2/3({n-1/24})}}} {4\sqrt{3}\lambda_n^2}$ 
\\  $+O(e^{\pi \sqrt{2/3({n-1/24})}}/
\lambda_n^3),
$  
 \\
$\lambda_n=\sqrt{n-1/24}$, 
\end{tabular}
& 
{}
\begin{tabular}{l}
Good bases,  \\
finite dimensions,  \\
many relations \\ and identities
\end{tabular}
&
{}
\begin{tabular}{l}
Values \\
of $L$-functions, \\ 
periods, \\
congruences, \dots
\end{tabular}
\end{tabular}
\noindent

Other examples:
Birch and Swinnerton-Dyer conjecture, \dots 
$L$-values attached to modular forms

\subsubsection*{Our data: three primitive cusp eigenforms}

{
\begin{align}\label{EQfj}
 f_j(z)=\sum_{n=1}^\infty a_{n,j}q^n\in \Sr_{k_j}(N_j, \psi_j),\  (j=1, 2,3)
\end{align}
}of weights  $k_1, k_2, k_3$, of 
 conductors $N_1, N_2, N_3$,
 and of nebentypus characters  $\psi_j \bmod N_j$, $N:=  \lcm(N_1,N_2,N_3)$.

Let  $p$ be a prime, $p\nmid N$.

We view $f_j\in \Qb[\![q]\!]\buildrel i_p\over\hookrightarrow{\mathbb C}_p[\![q]\!]$
 via a fixed embedding
$\Qb\buildrel i_p\over\hookrightarrow{\mathbb C}_p$, ${\mathbb C}_p={\widehat\Qb}_p$
 is Tate's field.
{} 


Let $\chi$ denote a {\sf{} variable} Dirichlet character ${}\bmod Np^v, v\ge 0$.

We view $k_j{}$ as a {\sf{} variable} weight 
in the {\sf {} weight space $X=X_{Np^v}
=\Hom_{cont}(Y, {\mathbb C}_p^*)$, $Y=(\Z/N\Z)^*\times\Z_p^*\ni(y_0, y_p)
$}.

{}{{}
The space $X$ is a {{} $p$-adic analytic space} first used in Serre's \cite{Se73}
{\em Formes modulaires et fonctions z\^eta $p$-adiques}.}
{}
 {\sf{}
Denote by  $(k, \chi)\in X$  the homomorphism $(y_0, y_p)\mapsto \chi(y_0)\chi(y_p\bmod p^v)y_p^k$.}
{}
We write simply $k_j$ for the couple $(k_j,\psi_j)\in X$.

\subsubsection*{
The purpose of this paper is to describe a  four variable
 $p$-adic
$L$ function}
\noindent attached to Garrett's triple product of three Coleman's families
{
$$
k_j{}\mapsto \left \{f_{j,k_j{}}= \sum_{n=1}^\infty a_{n,j}(k{}_j)
q^n\right \} 
$$  
}of cusp eigenforms of three constant slopes $\sigma_j=\ord_p(\alpha_{p, j}^{(1)}(k_j{}))\ge 0$
 where {{} $\alpha_{p, j}^{(1)}(k_j)$, $ \alpha_{p,j}^{(2)}(k_j)$}  
are  the {\sf{} Satake parameters} 
given as inverse roots of  
the Hecke $p$--polynomial
$
 1-a_{ p,j}X-\psi_j(p)p^{k_j-1}X^2 = (1-\al_{ p,j}^{(1)}(p)X)
(1-\al_{ p,j}^{(2)}(p)X). \ \ 
 $

We assume that {{}$\ord_p(\alpha_{p, j}^{(1)}(k_j{}))\le\ord_p(\alpha_{p, j}^{(2)}(k_j{}))$}. 


{This extends a previous result}:
(see
\cite{PaTV},  
 Invent. Math. v. 154, N3 (2003)) 
where a {\sf two variable $p$-adic $L$-function} was constructed
interpolating on all $k{}$ a function  $(k{},s)\mapsto L^*(f_{k{}}, s, \chi)$ $(s=1, \cdots, k{}-1)$
for such a family. 


We use the theory of $p$-adic integration 
 with values 
 in spaces of
 {\sf{} nearly holomorphic  modular forms}
 (in the sense of Shimura, see \cite{ShiAr}).



\subsubsection*
{  A
family of slope $\si>0$ of cusp eigenforms $f_{k{}}$ of weight $k{}\ge 2$:}

\begin{tabular}{ll}
\begin{tabular}{l}{}
\begin{tabular}{l}
$\ds{}
k{}\mapsto f_{k{}}= \sum_{n=1}^\infty a_n(k{})q^n$\\
$\in \Qb[\![ q]\!] \subset {\mathbb C}_p[\![ q]\!] $ \\ 
\end{tabular}\\ {}
\begin{tabular}{l}
A model example  \\ 
of a $p$-adic
family \\ 
({\sf  not cusp and $\si=0$}): 
\end{tabular}
\\ {}
\begin{tabular}{l}
{\sf Eisenstein series}\\
$\ds{}
a_n(k{}) = \sum_{d|n}d^{k{}-1}, 
f_{k{}}=E_{k{}}
$

\end{tabular}
\end{tabular}
&
\hskip-0.5cm
\begin{tabular}{|l}

\begin{tabular}{l}
1)
the Fourier coefficients $a_n(k{})$ of $f_{k{}}$ \\
and one of 
the Satake 
\\ 
$p$-parameters $\alpha(k{}):=\alpha_p^{(1)}(k{})$  
\\
  are given by 
certain 
 $p$-adic analytic 
\\ 
{ functions  $k{}\mapsto a_n(k{})$ 
for $(n,p)=1$} 
\end{tabular}\\ {}
\begin{tabular}{l}\\
2) 
the slope is {\sf constant 
and positive}: \\  
{} $\ord(\alpha(k))=\si > 0$  \\
\end{tabular}
\end{tabular}
\end{tabular}


{\sf  The existence of  
families of slope $\si>0$  was established
 in \cite{CoPB}}

\begin{tabular}{ll}
\hskip-0.5cm
\begin{tabular}{l}
R.Coleman gave an example with \\ $p=7$, $f=\De$, $k=12$
\\ 
$a_7=\tau(7)=-7\cdot 2392, \si=1$. 
\end{tabular}
&{}
\hskip-0.5cm
\begin{tabular}{|l}
A program in PARI for computing \\
such families is 
contained in 
\cite{CST98}
\\
(see also the Web-page of W.Stein,
\\
{\tt http://modular.fas.harvard.edu/}
)
\end{tabular}
\end{tabular}
\subsubsection*{Coleman proved that :}

{}
 \begin{tabular}{ll}
 \begin{tabular}{l}
\noindent$\bullet$
The operator $U$
acts as \\
a completely continuous operator
\\ on each $\Ar$-submodule 
${\cal M}^\dagger (N p^{v}; \Ar)$\\ $\subset \Ar[\![ q]\!]$
 (i.e. $U$ is a limit \\ of
 finite-dimensional operators) 
 \end{tabular} & 
 \begin{tabular}{l} 
$\Longrightarrow $ there exists \\
the {\sf{} Fredholm determinant} \\ $P_U(T)$\\ $=
\det(Id-T\cdot U)\in \Ar[\![ T]\!]$
\\ \\ 
\end{tabular}
 \end{tabular}


 \begin{tabular}{ll}
 \begin{tabular}{l}
\noindent$\bullet$
there is a version \\ of the {\sf{} Riesz theory}:\\
for any inverse root  $\al\in \Ar^*$ \\ of $P_U(T)$
there exists \\ an eigenfunction $g$, $Ug=\al g$
 \end{tabular} & \qquad	
 \begin{tabular}{l} such that 
$ev_{k}(g)\in {{\mathbb C}_p}[\![ q]\!]$ \\ are classical cusp eigenforms \\ for all 
$k$  
in a neigbourhood
\\ ${\cal B}\subset
 {X}$
(see in \cite{CoPB})
\end{tabular}
 \end{tabular}



\section{Generalities on triple products}
\label{ptr0-1}

The triple product  with a Dirichlet character $\chi$ is defined 
 as the following  complex  $L$-function  
({\sf{} an Euler product of degree eight}):
{{}
  \begin{align}\label{ptr0-1.1} &
 L(f_1\otimes f_2\otimes f_3, s, \chi) = 
\prod_{p\nmid N}
 L((f_1\otimes f_2\otimes f_3)_p, \chi(p) p^{-s}), 
\end{align}}
 \begin{align} \label{ptr1.2} &{}
\mbox{ where } 
{}
 L((f_1\otimes f_2\otimes f_3)_p, X)^{-1}  = 
\\ 
& {} {{}
\det \left(1_8 -  X\mat \al^{(1)}_{ p,1}, 0, 0, \al^{(2)}_{ p,1}, \otimes
 \mat \al^{(1)}_{ p,2}, 0, 0, \al^{(2)}_{ p,2}, \otimes\mat 
\al^{(1)}_{ p,3}, 0, 0, 
\al^{(2)}_{ p,3},
 \right) } \nonumber\\ 
&{}
  = \prod_{\eta}(1-\al_{ p,1}^{(\eta(1))}\al_{ p,2}^{(\eta(2))}
\al_{ p,3}^{(\eta(3))}X) \ \ \
\nonumber\\ 
&{}
  =(1-\al_{ p,1}^{(1)}\al_{ p,2}^{(1)}
\al_{ p,3}^{(1)}X)(1-\al_{ p,1}^{(1)}\al_{ p,2}^{(1)}
\al_{ p,3}^{(2)}X)\cds (1-\al_{ p,1}^{(2)}\al_{ p,2}^{(2)}
\al_{ p,3}^{(2)} X), 
\nonumber
 \end{align}

product taken over all  8 maps $\eta: \{1, 2, 3 \}\to \{1, 2\}$.
\subsubsection*{The Satake parameters and  Hecke $p$--polynomials of forms  $f_j$:}
Here the Satake parameters 
$\al_{ p,j}^{(1)}$, $\al_{ p,j}^{(2)}$ 
are given as inverse roots of  
the Hecke $p$--polynomials

 $$
 1-a_{ p,j}X-\psi_j(p)p^{k_j-1}X^2 = (1-\al_{ p,j}^{(1)}(p)X)
(1-\al_{ p,j}^{(2)}(p)X),\ \ 
 $$

 We always assume that 
the weights are ``balanced'': 

 \begin{eqnarray}\label{ptr0-1.3}
k_1 \ge k_2\ge k_3\ge 2,  \mbox{ and } k_1 \le k_2+ k_3-2 
\end{eqnarray}



\subsubsection*{Critical values and functional equation}
We use the corresponding normalized $L$ function
 (see \cite{De79}, \cite{Co}, \cite{Co-PeRi}),
 which 
 has the form: 

\begin{align}\label{ptr6.4} & 
\La (f_1\otimes f_2\otimes f_3, s, \chi) = 
\\ &
\Ga_{\mathbb C}(s)
\Ga_{\mathbb C}(s-k_3+1)\Ga_{\mathbb C}(s-k_2+1)\Ga_{\mathbb C}(s-k_1+1)
 L(f_1\otimes f_2\otimes f_3, s, \chi),\nonumber
\end{align}

where {{}$\Gamma _{{{\mathbb C}}}(s) = 2(2\pi )^{-s}\Gamma (s)$}.

 The Gamma-factor 
determines the {\sf{} critical values}  
$s=k_1, \cdots, k_2+ k_3-2$ of $\La(s)$, which we explicitely
evaluate (like in the classical formula $\zeta(2)=\ds \frac {\pi^2} 6$). 

{\sf{} A  functional equation} of $\La(s)$
has the form:
 $$
s\mapsto k_1+k_2+k_3-2-s.
$$

According to H.Hida \cite{Hi86} and R.Coleman \cite{CoPB},
 one can include each $f_j$ $(j=1, 2, 3)$
 (under suitable  assumptions on $p$ and on $f_j$)  into a
 {{} $p$-adic analytic family}
$$
\textbf{f}_j : k_j{}\mapsto \{f_{j,k_j{}}= \sum_{n=1}^\infty a_{n}(f_{j, k_j{}})q^n\}
$$ 
 of cusp eigenforms $f_{j,k_j{}}$ of weights $k_j{}$ 
 in such a way that $f_{j,k_j}=f_j$,
  and that all their Fourier coefficients $a_n(f_{j, k_j{}})$
  are given by 
 certain $p$-adic analytic functions
$k_j{}\mapsto a_{n, j}(k_j{})$.

\section{Statement of the problem}

\subsubsection*{Given three $p$-adic analytic families 
$\textbf{f}_j$ of slope $\si_j\ge0$, 
to construct a four-variable $p$-adic $L$-function
attached to  Garrett's triple product
of these families.}
\noindent
We show that this function 
 interpolates the special values
$$ 
(s, k_1{}, k_2{}, k_2{})\longmapsto
\Lambda (f_{1,k_1{}}\otimes f_{2,k_2{}}\otimes f_{3,k_3{}}, s, \chi)
$$
{}
at critical points $s=k_1{}, \cdots, k_2{}+k_3{}-2$
for balanced weights $k_1{} \le k_2{}+k_3{}-2$;
 we prove that these values are 
algebraic numbers after dividing out certain ``periods''.

{}
However the construction uses directly modular forms, 
and not the $L$-values in question, 
and a comparison of special values of two functions is
done {\sf{} after the construction}. 
{}






Consider the product of the Satake parmeters 
$$
\lambda_p=\al^{(1)}_{ p,1}\al^{(1)}_{ p,2}\al^{(1)}_{ p,3}=\lambda_p(k_1{}, k_2{}, k_3{}) 
$$ 
{}
We  assume that
 {{}
$
\ord_p \al^{(1)}_{ p,j}\le \ord_p\al^{(2)}_{ p,j}$, 
}
 and that the slope 
{}
$
\sigma= \ord_p(\lambda_p(k_1{}, k_2{}, k_3{}))
$
{}
 is {\sf{} constant and positive } for all triplets $(k_1{}, k_2{}, k_3{})$ 
 in a 
$p$-adic neighbourhood $\Br\subset X^3$
 of the fixed triplet of weights $(k_1, k_2, k_3)$.

\subsubsection*{  Our method includes:}


$\bullet$ a version 
 of 
{\sf {} Garrett's integral representation} 
for the triple $L$-functions of the form: 
{}
for $r=0, \cdots, k{}_2+k{}_3-k{}_1-2$, \\
$\ds \La (f_{1,k{}_1}\otimes f_{2,k{}_2}\otimes f_{3,k{}_3}, k{}_2+k{}_3-r, \chi)
 = $
\\
$\ds
\mathop{\int\int
\int}_{\left(\Gamma_0(N^2p^{2v})\backslash {\HH}\right)^3}
 \overline{{\tilde f_{1,k{}_1}}(z_1){\tilde f_{2,k{}_2}}(z_2){\tilde f_{3,k{}_3}}(z_3)}
{\cal E}({z}_1,{z}_2,{z}_3;-r, \chi)\prod_j ({dx_jdy_j\over y_j^2})
$
\\ 
where {{}\sf $\tilde f_{j,k_j}=:f_{j,k_j}^{\,0}$} is an eigenfunction of 
$U_p^*$ in ${\cal M}_{k_j}(Np, \psi_j)$,

{{} $f_{j,k_j, 0}$} is the corresponding eigenfunction of $U_p$,
\\

{{}${\cal E}({z}_1,{z}_2,{z}_3;-r, \chi)\in\scrM_T(N^2p^{2v})$

$=\scrM_{k_1{}}(N^2p^{2v}, \psi_1)\otimes \scrM_{k_2{}}(N^2p^{2v}, \psi_2)\otimes 
\scrM_{k_3{}}(N^2p^{2v}, \psi_3)
$ }\\ 
is the 
triple modular form of triple weight $(k{}_1,k{}_2, k{}_3)$, and of 
fixed triple Nebentupus character $(\psi_1, \psi_2, \psi_3)$, 
obtained from a {}
nearly  holomorphic Siegel-Eisenstein series
$F_{\chi, r}=
G^{\star}({z}, -r; k{}, {(Np^v)}^2,\bpsi),
$
of degree 3, of weight $k{}=k_2+k_3-k_1$, 
and the Nebentypus character 
$\bpsi=\chi^2\psi_1\psi_2\overline\psi_3$

\subsubsection*{
We obtain ${\cal E}({z}_1,{z}_2,{z}_3;-r, \chi)$ from a Siegel-Eisenstein series}
by applying to $F_{\chi, r}$ {{} Boecherer's higher twist} (see  (\ref{Fchichi}))
and  {{} Ibukiyama's differential operator} (see (\ref{IbuOp})).

{
These operations act explicitely on the Fourier expansions.}

Then one uses:

{}
$\bullet$
The   {\sf {} theory of $p$-adic integration} 
 with values 
in Serre's type  $\Ar$-modules $\Mr_T(\Ar)$ 
 of {{} triple arithmetical nearly holomorphic  modular forms}
 over $p$-adic Banach algebras $\Ar$.
Explicit Fourier coefficients $a_{\chi,r}(R,{\cal T})\in \Qb[R,T]$
of 
${\cal E}(-r, \chi)$ are given by special polynomials 
of matricies ${\cal T}=(t_{ij})$, $R=(R_{ij})$
and of $\chi(\beta)\beta^r$ (with $\beta\in \Z_p^*\cap\Q$) i.e. the coefficients of $a_{\chi,r}$
 by some elementry $p$-adic measures $\ds\int_Y\chi y^r\dd \mu_{\cal T}\in {\cal A}$.
Here  $\Ar=\Ar(\Br)$ is a certain $p$-adic Banach algebra
of {\sf{}  functions on an open analytic subspace 
$\Br=\Br_1\times\Br_2\times\Br_3\subset X^3$ in the product of three copies
 of the weight space} $X=\Hom_{cont}(Y,{\mathbb C}_p^*)$.

{}{{}
These  measures  on the group $Y= (\Z/N\Z)^* \times\Z_p^*$ produce the coefficients of $a_{\chi,r}$
of ${\cal E}(-r, \chi)$ of $\Mr_T(\Ar)$
for all $p$-adic weights $x\in X$, 
given by
$\ds\int_Yx(y)\dd \mu_{\cal T}\in {\cal A}$  ({\it an interpolation from $x=\chi y_p^r$ to all $x\in X$}).
} 

\subsubsection*{
$\bullet$
The {\sf {} spectral theory of triple Atkin's operator $U=U_ {p,T}$}
}

allows to 
evaluate the integral using 
at each weight $(k_1, k_2, k_3)$ the equality
${\ang{\underline f^{\,0}, {\cal E}(- r, \chi)},}=\ang{\underline f^{\,0}, \pi_\lambda({\cal E}(- r, \chi))},$ with the
projection $\pi_\lambda$ of $\Mr_T(\Ar)$ to the $\lambda$-part ${\Mr_T}(\Ar)^\lambda$,  defined by :
{
$
\Ker \pi_{\lambda}:=\bigcap_{n\ge 1}\Im 
(U_T-{\lambda} I)^n, 
$
$
\Im \pi_{\lambda}:=\bigcup_{n\ge 1}\Ker
(U_T-{\lambda} I)^n.
$
}

We prove that 
$U$ is a completely continuos $\Ar$-linear operator
on a certain Coleman's submodule 
${\cal M}({\cal A})^\dagger$ of Serre's type module ${\cal M}({\cal A})$.
Then  the projection 
$\pi_\lambda$ exists (on this submodule)
due to general results of Serre and Coleman, see \cite{CoPB}, \cite{SePB}.
{}

{{}
We show that there exists an element 
$\tilde {\cal E}(- r, \chi)\in{\cal M}({\cal A})^\dagger$ such that
at each weight $(k_1, k_2, k_3)$ the equality holds:\\
${\ang{\underline f^{\,0}, {\cal E}(- r, \chi)},}=\ang{\underline f^{\,0}, \pi_\lambda(\tilde{\cal E}(- r, \chi))},$, 
and the product can be expressed through certain coefficients 
the series $\tilde{\cal E}(- r, \chi)$ which are the same as those of ${\cal E}(- r, \chi)$.
}


\subsubsection*{$\bullet$  Key point: modular admissible measures}

Let us write for simplicity: {{}
${\cal E}(- r, \chi)$ for $\tilde{\cal E}(- r, \chi)$

${\cal M}_T({\cal A})$ instead of
${\cal M}_T({\cal A})^\dagger$}  (Coleman's submodule)

One  defines    {\sf {} admissible $p$-adic measures} $\tilde\Phi^\lambda
$
with values 
in Banach  $\Ar$-modules $\Mr_T^\lambda(\Ar)$ which are locally free of finite rank, 
using the {\sf{} test functions:
$\ds\int_Y\chi y_p^r\tilde\Phi^\lambda=
\pi_\lambda({\cal E}(-r, \chi))$}.

Consider
the {\sf{{} evaluation maps
$ev_{\textbf{s}}:\Ar\to {\mathbb C}_p$  for any $p$-adic triple weights
$\textbf {s}=(s_1, s_2, s_3)\in \Br$}}.


\noindent$\bullet$ {{} Passage from values in modular forms to scalar values}:
apply an algebraic
$\Ar$-linear  form $\Mr_T^\lambda(\Ar)\buildrel \ell_T\over\to\Ar$
 to the constructed measure $\tilde\Phi^\lambda$ (in modular forms), 
and the  {\sf{} evaluation maps}
$\Ar\buildrel{ev_{\textbf{s}}}\over\to {\mathbb C}_p$  for any $p$-adic triple weights
$\textbf {s}\in X^3$.

{}
The linear form $\ell_T$ is an algebraic version of the Petersson product
(a geometric meaning of $\ell_T$: the first coordinate in an (orthogonal) $\Ar$-basis of 
eigenfunctions of all Hecke operators $T_q$ for $q\nmid Np$, 
with the first basis element $\textbf {f}_0\in \scrM^\lambda(\scrA)$). 


\subsubsection*{Using the evaluation map and the Mellin transform}

We obtain  the { measure  $\mu=\ell_T( \tilde\Phi^\lambda)$  with values in  
$\Ar$ on the profinite
group $Y$}.{}

 \noindent$\bullet$
Construct an analytic function  
$
\Lr_\mu:X\to \Ar=\Ar(\Br)  
$
 as the {\sf{} $p$-adic Mellin transform}
$\ds
\Lr_\mu (x)=\int_{Y} x(y)\, d\mu(y)\in\Ar$, $x\in X$. 


{}
 \noindent
$\bullet$
{{} Solution}: the  function in question
$\Lr_\mu (x, {\textbf{s}})$ is given 
by evaluation of $\Lr_\mu (x)$ at ${\textbf{s}}=(s_1, s_2, s_3)\in \Br$:
this is  a 
$p$-adic analytic function 
 in four variables 
$$
{}
(x, {\textbf {s}})
\in X\times \Br_1\times\Br_2\times\Br_3\subset X\times X
\times X
\times X
 $$
$
\Lr_{\mu} (x, {\textbf s}) := 
ev_{\textbf {s}}(\Lr_{\mu} (x)) \ \ \
 (x \in X, \ 
{\textbf {s}}
\in  \Br_1\times\Br_2\times\Br_3, 
\  \Lr_{\mu} (x)\in \Ar).
$

\subsubsection*{Final step: comparison between ${\mathbb C}$ and ${\mathbb C}_p$}

 \noindent$\bullet$
We check an  equality relating the values
of the constructed analytic function
{\sf
$
\Lr_{\mu} (x, {\textbf s}) 
$
at the  arithmetical characters} \\ $ x=y_p^r\chi\in X$, 
and {\sf{} at triple weights ${\textbf s}= ({k{}_1, k{}_2, k{}_3})\in\Br$}, 
with the normalized critical special values
{}
$$
L^*({f_{1,k_1{}}}\otimes {f_{2,k_2{}}}
\otimes {f_{3, k_3{}}}, k_2{}+k_3{}-2-r, \chi)  \ \ (r=0, \cdots, k_2{}+k_3{}-k{}_1-2), 
$$

 for  certain Dirichlet characters
$\chi\bmod Np^v, v\ge 1$.
These are  {\sf{} algebraic numbers},  
embedded into  {{}${\mathbb C}_p=\widehat{\overline\Q}_p$ (the Tate field 
of $p$-adic numbers).
}
The normalisation of $L^*$ includes at the same time
Gauss sums, Petersson scalar products, 
powers of $\pi$,  the product $\lambda_p(k{}_1, k{}_2, k{}_3)$, 
and a certain finite Euler product.


\section{Arithmetical  nearly holomorphic modular forms} 
\subsubsection*{Arithmetical  nearly holomorphic modular forms (the elliptic case)}
Let $\Ar$ be a { commutative ring} (a subring of ${\mathbb C}$ or ${\mathbb C}_p$)


{\sf
Arithmetical nearly holomorphic  modular forms} (in the sense of Shimura, \cite{ShiAr}
  are certain formal series

\begin{align*}&
g=\sum\limits_{n=0}^\infty a(n; R)q^{n}
 \in 
{\Ar}[\![q]\!][R], \mbox{ with the property }
\end{align*}
{}
 that for ${\Ar}={\mathbb C}$, $z=x+iy\in \dbH, \ R=(4\pi y)^{-1}$, 
  the series converges to a $\scrC^\infty$-modular form on 
$\dbH$ of a given
 weight $k$ and Dirichlet character $\psi$.
The coefficients $a(n; R)$
 are polynomials in ${\Ar}[R]$.
If $\deg_R a(n; R)\le r$ for all $n$, we call $g$ 
{\sf{} nearly holomorphic of type }$r$
(it is annihilated by $(\frac \partial{\partial \overline z})^{r+1}$, 
see \cite{ShiAr}).

We use the notation 
{
$\Mr_{k, r}(N, \psi, \Ar)$ or $\tilde \Mr(N, \psi, \Ar)$
}
 for $\Ar$-modules of such forms 
(In our constructions the weight $k$ varies).

A {known example} (see the introduction to \cite{ShiAr})
is given by
the series
{
\begin{align*}&
-12R+E_2:=
-12R+1
-24\sum_{n=1}^\infty \si_{1}(n) q^n
\\ &
=\frac 3{\pi^{2}}
\lim_{s\to 0}
y^s\sum_{m_1, m_2 \in \Z}\!\!\!\!{}^{\displaystyle {}'} 
(m_1+ m_2 z)^{-2}|m_1+ m_2 z|^{-2s}, (R=(4\pi y)^{-1})
\end{align*}}where $\si_{1} (n) = \sum_{d|n}d^{}$.

The action of the {\sf{} Shimura differential operator}
$$
\delta_k: \Mr_{k, r}(N, \psi, \Ar)\to\Mr_{k+2, r+1}(N, \psi, \Ar), 
$$
is given over ${\mathbb C}$ by  $\ds \delta_k(f)=(\frac 1{2\pi i}\frac\partial{\partial z}-\frac k {4\pi y})f$.

This operator is a correction  of the  {\sf{} Ramanujan operator} 
\begin{align*}&
\theta (\sum_{n=0}^\infty a_nq^{n})=
\sum_{n=1}^\infty na_nq^{n}
=
 \frac 1{2\pi i}\frac{\partial }{\partial z}(\sum_{n=0}^\infty a_n
q^{n})=q\frac {\partial }{\partial q}(\sum_{n=0}^\infty a_n
q^{n}), 
\end{align*}
which does not preserve the modularity.
For example  
 $\theta\Delta=E_2\Delta$, where $E_2$ is a 
{\sf{} quasimodular} form
(in the sense of Kaneko and  Zagier, see
\cite{Ka-Za}).

Notice that $\delta_k f= (\theta -k R)f$, and that 
 {\sf{} Serre's operator} 
$
f\mapsto \theta f -\frac k {12} E_2f$
takes  $\Mr_k$ to  $\Mr_{k+2}$.

Note that that the {\sf{} arithmetical twist operator} 
\begin{align*}&
\theta_\chi (\sum_{n=0}^\infty a_nq^{n})=
\sum_{n=1}^\infty \chi(n)a_nq^{n}
\end{align*}
is a natural  analog of  
the  Ramanujan operator.

\subsubsection*{Triple arithmetical modular forms}
{Let $\Ar$ be a commutative ring.}
The tensor product over $\Ar$  
$$
\Mr_{{\textbf k}, r, T}(N, \psi, \Ar):=
\Mr_{k_1, r}(N, \psi_1, \Ar)\otimes\Mr_{k_2, r}(N, \psi_2, \Ar)
\otimes\Mr_{k_3, r}(N, \psi_3, \Ar)
$$
consists of {\sf {} triple arithmetical modular forms} as 
 certain formal series of the form 
{
\begin{align*}&
g=\sum\limits_{n_1, n_2, n_3=0}^\infty a(n_1, n_2, n_3; R_1, R_2, R_3)q_1^{n_1}q_2^{n_2}q_3^{n_3}
\\ & \in 
{\Ar}[\![q_1, q_2, q_3]\!][R_1, R_2, R_3], \mbox{ where }z_j=x_j+iy_j\in \dbH, \ R_j=(4\pi y_j)^{-1}, 
\end{align*}}{}
with the property that for ${\Ar}={\mathbb C}$, 
  the series converges to a $\scrC^\infty$-modular form on $\dbH^3$ of a given
 weight $(k_1,  k_2, k_3)$ and character $(\psi_1, \psi_2, \psi_3)$, $j=1, 2, 3$.
The coefficients $a(n_1, n_2, n_3; R_1, R_2, R_3)$
 are polynomials in ${\Ar}[R_1, R_2, R_3]$.
Examples of such modular forms come
 from the restriction to the diagonal 
of Siegel modular forms of degree 3.

\section{Siegel-Eisenstein  series}

\subsubsection*{
Siegel
modular groups}
{}
Let $\ds J_{2m}=\begin{pmatrix}
0_m &-1_m\cr
1_m &0_m\end{pmatrix}$.{}
The symplectic group
{}
$$
\Sp_{m}(\dbR)=\left\{g\in\GL_{2m}(\dbR)\vert
{}^t g\cdot J_{2m}g=J_{2m}\right\}, 
$$
{}
acts on the Siegel upper half plane 
$$
\HH_m=\left\{z={}^t z\in M_m(\dbC)\vert
\im z>0\right\}
$$
by $g(z)=(az+b)(cz+d)^{-1}$, 
where we use the
bloc notation $g=\left({a\,\,b\atop
c\,\,d}\right)\in \Sp_{2m}(\R)$. {} 
We use the 
{{}
congruence subgroup
$\Gamma_0^m(N)=\{\ga\in \Sp_{m}(\Z)\ |\ \ga\equiv\mat  *,*,0,*, \bmod N
\}\subset \Sp_{m}(\Z)$.
}

\subsubsection*{A Siegel modular form} 
{}
$f\in
\Mr_k(\Gamma_0^m(N),\chi)$ of degree $m>1$, 
weight $k$ and a Dirichlet chracter 
$\chi\bmod N$ {}
is a holomorphic
function $f\colon\, \HH_m\to\dbC$ such that
for every $\gamma=\left({a\,\,b\atop
c\,\,d}\right)\in\Gamma_0^m(N)$ one has {}
{{}
$$
f(\gamma(z))=\chi(\det\,d)\det(cz+d)^k
f(z).
$$
}

The {\sf  Fourier expansion}  of  $f$ uses
the symbol
$q^{\scrT}=\exp(2\pi i {\rm tr}({\scrT}z)) $ \\
$
= \prod_{i=1}^mq_{ii}^{{\scrT}_{ii}}\prod_{i<j}q_{ij}^{2{\scrT}_{ij}}\subset 
\dbC[\![q_{11}, \dots , q_{mm}]\!][q_{ij}, \ q_{ij}^{-1}]_{i,j=1, \cdots, m},
$, \\
$q_{ij}= \exp(2\pi (\sqrt {-1} z_{i,j}))$, 
{} and
$\scrT$ in the semi-group
 {} $
B_m =\{{\scrT}={}^t{\scrT}\ge 0\vert {\scrT}\hbox{
half-integral}\}: 
$\\
$
f(z) =\sum\limits_{{\scrT}\in
B_m}a({\scrT})q^{\scrT}\in\dbC[\![q^{B_m}]\!]
  (\mbox{a formal }
q\mbox{-expansion }\in {\mathbb C}[\![q^{B_m}]\!]),
${}\\

\subsubsection*{Siegel-Eisenstein series}

\begin{exam}[ {
 \cite{Nag2}, p.408}] 

\begin{align*}
E^{(2)}_4(z)=&1+240q_{11}+240q_{22}+2160q^2_{11}+ (240q_{12}^{-2}
+13440q_{12}^{-1} \cr
 &  +30240+ 13440q_{12}+ 240q_{12}^{2})q_{11}q_{22}+2160q^2_{22}+\dots 
\\  
E^{(2)}_6(z)= &1-504q_{11}-504q_{22}-16632q^2_{11}+ (-540q_{12}^{-2}
+44352q_{12}^{-1} \cr 
 &  +166320+ 44352 q_{12}-504q_{12}^{2})q_{11}q_{22}-16632q^2_{22}+\dots .
\end{align*}
\end{exam}


\subsection*{Arithmetical nearly holomorphic Siegel modular 
forms}
\subsubsection*{Arithmetical  Siegel modular 
forms}

Consider  
a  commutative ring $\Ar$, 
the  formal variables $q=(q_{i,j})_{i,j=1, \dots, m}$, $R=(R_{i,j})_{i,j=1, \dots, m}$,  
and the ring of {\em formal Fourier series}
{{}
\begin{align}\label{ArFS}
\Ar [\![q^{B_m}]\!][R_{i,j}]=\left\{f =\sum\limits_{{\scrT}\in
B_m}a({\scrT}, R)q^{\scrT}\ \Big|\ a({\scrT}, R)\in \Ar[R_{i,j}]\right\}
\end{align} 
}
(over the complex numbers this notation corresponds to
 $q^{\scrT}=\exp(2\pi i {\rm tr}({\scrT} {z}))$,  
$R=(4\pi{\rm Im}({z}))^{-1}$).

{}
The formal Fourier expansion of a nearly holomorphic 
Siegel modular 
form $f$ with coefficients in $\Ar$ 
is a certain  element of $\Ar [\![q^{B_m}]\!][R_{i,j}]$.
We call $f$ {\sf{} arithmetical} in the sense of Shimura
\cite{ShiAr}, if $\Ar=\Qb$. 

\subsection{Algebraic differential operators of Maass and Shimura}\label{sd2}
\subsubsection*{Maass
differential operator}
Let us consider the {\sf{} Maass
differential operator} (see \cite{Maa71})
$\De_m$  of degree $m$, acting on complex
${\mathcal C}^{\infty}$-functions on $\dbH_m$ by:
\begin{align}\label{sde2.1}&
    \De_m    = \det(\tpa_{ij}), \hskip1cm
    \tpa_{ij} = 2^{-1}(1+\de_{ij}) \pa/\pa_{ij}, 
\end{align}{}
its algebraic version is the {\sf{} Ramanujan operator of degree $m$}:  
\begin{align}\label{sde2.1}&
\Theta_m:=\det(\frac 1{2\pi i} \tpa_{ij})=\det(\theta_{ij})=\frac 1{(2\pi i)^m}\De_m, 
\end{align}{}
where 
$\Theta_m(q^\scrT)=\det(\scrT)q^\scrT.$

\subsubsection*{Shimura differential operator}


The
{\sf{} Shimura differential operator} (see \cite{Shi76, ShiAr}):
$$  \de_k f(z)  =\det(R)^{k+1-\vkp} \Theta_m \left[
                  \det(R)^{\vkp-1-k} f \right], \mbox{ where
}
 R=(4\pi y)^{-1}, 
$$
acts on arithmetic nearly holomorphic Siegel modular forms, and
the composition is defined
{}
\begin{align}\label{dek}&
\de_k^{(r)} = \de_{k+2r-2} \circ \cdots \circ \de_k :
\Mt^m_k(N, \psi;\Qb) \rightarrow \Mt^m_{k+2rm}(N, \psi;\Qb), 
\\ \nonumber & 
\mbox{ where }
\\ \nonumber & 
\de_k f(z)  =  
         \Kpi^m \det(y)^{-1} \det(z-\bar{z})^{\vkp-k} \De_m \left[
                  \det(z-\bar{z})^{k-\vkp+1} f(z) \right].
\end{align}

\subsubsection*{Universal polynomials $Q(R,{\scrT}; k, r)$}

Let 
$\ds f = \sum_{{\scrT}\in B_m} c({\scrT})q^{\scrT}  \in \Mr_k^m(N, \psi)$
be a formal holomorphic Fourier  expansion.
One  shows   that 
$\de_k^{(r)} f$
is given by
{}
$$
\de_k^{(r)}f =
\sum_{{\scrT} \in B_m}
Q(R,{\scrT}; k, r) c({\scrT}) 
q^{\scrT}. 
$$
\subsubsection*{Universal polynomials (continued)}

Here we use a universal polynomial (\ref{QRT})
which can be defined for all $k\in {\mathbb C}$, and 
it expresses the action of the Shimura operator 
on the exponential (of degree $m$):
$$
\de_k^{(r)}(q^{\scrT})=Q(R,{\cal T}; k, r)
q^{\scrT}.
$$
If $m=1$, $r$ arbitrary (see \cite{Shi76}),
{}
$\ds
\delta_k^{(r)}=\sum_{j=0}^r(-1)^{r-j}{r\choose j}\frac {\Gamma(k+r)}{\Gamma(k+j)}
R^{r-j}\theta^j$, \\ $\ds Q(R,n; k, r)=\sum_{j=0}^r(-1)^{r-j}{r\choose j}\frac {\Gamma(k+r)}{\Gamma(k+j)}
R^{r-j}n^j.
$
{}\

\subsubsection*{Universal polynomials (continued)}

If $r=1$, $m$ arbitrary,  one has (see \cite{Maa71}):
$$
{}
\de_k f(z)  = \sum_{\scrT \in B_m} c(\scrT) 
  \sum_{l=0}^m (-1)^{m-l} 
c_{m-l}(k+1-\vkp) \tr \left( {}^t \rho_{m-l}(R )
                           \cdot \rho_{l}^{\star}(\scrT)\right )q^\scrT 
$$
where $R=(4\pi y)^{-1}=(R_{i,j})\in\M_m(\R)$, 
$\di{ c_m(\al) = \fragam{\al+\ka}{\al+\ka-1} }$, 
{}
$\ds
\Gamma_m(s)=\pi^{m(m-1)/4}
\prod\limits_{j=0}^{m-1}\Gamma(s-(j/2))$.

Here we use
{
 the natural representation
$
       \rho_r    :    \GL_m({\mathbb C})    \longrightarrow    \GL(\wedge^r{\mathbb C}^m) 
$
 ($0 \leq r \leq m$)
of the group $\GL_m(\dbC)$ on the vector space $\Lb^r\dbC^m$.
}
Thus $\rho_r(z)$ is a matrix of size ${m \choose r}\times {m \choose r}$ composed of the
{\sf{} subdeterminants of $z$ of degree $r$}.
Put
${}
\rho_r^{\star}(z) = \det(z) \rho_{m-r}({}^t\! z)^{-1}{}.
$

Then the representations $\rho_r$ and $\rho_r^{\star}$ turn
 out to be polynomial representations.

In general 
 (see \cite{CourPa}, 
Theorem 3.14) one has:
{}
\begin{align}\label{QRT}&
 Q(R,{\cal T})=Q(R,{\cal T}; k, r)
\\ &\nonumber \qquad\qquad =\sum_{t=0}^r{r\choose t} 
\det({\cal T})^{r-t}  \sum_{|L| \leq mt-t} 
R_L({\kappa} -k-r) Q_L(R,{\cal T}),
\\ &\nonumber
Q_L(R,{\cal T})=\tr \left( {}^t\! \rho_{m-l_1}(R )  \rho_{l_1}^{\star}
({\cal T} )\right) 
\cdot{\dots}\cdot \tr \left( {}^t\! \rho_{m-l_t}(R )  \rho_{l_t}^{\star}
({\cal T} )\right)). 
\end{align}
{}
In (\ref{QRT}), 
  $L$ goes over all the multi-indices
  $0 \leq l_1 \leq \cdots \leq l_t \leq m$, such that
  $|L| = l_1 + \cdots + l_t \leq mt-t$, and
  $R_L(\beta) \in \ZZ[1/2][\beta]$ in (\ref{QRT}) are {\sf{}
polynomials 
in $\beta$ of degree
  $(mt-|L|)$}
(used with $\beta={\kappa} -k-r$).

Note the {\sf {} differentiation  rule of degree $m$}  (see \cite{Shi83}, p.466):\\
$\ds
    \De(fg)    =     \sum_{r=0}^m \tr \left( {}^t\! \rho_r(\tpa/\pa z) f
                           \cdot \rho_{m-r}^{\star}(\tpa/\pa z) g \right).                    
$



\begin{exam}[Siegel-Eisenstein series of odd degree and higher level] 
{}
  \begin{align} \label{eqG*}&
G^*({z}, s; k, {\bpsi} , N) \\ \nonumber &=
\det(y)^s\sum_{c,d}{\bpsi}(
\det c)\det(c{z}+d)^{-k}
|\det(c{z}+d)|^{-2s}\
\cdot \\ &
\cdot  \tilde{\Gamma}(k,s) L_N(k+2s,{\bpsi}) 
                                    \left(
                                  \prod_{i=1}^{[m/2]} L_N(2k+4s-2i,{\bpsi}^2)
                                  \right)    \nonumber
, \mbox{ where 
}
  \end{align}
$(c,d)$ runs over all ``non-associated coprime symmetric pairs''
with $\det(c)$ coprime to $N$, 
 $\kp= (m+1)/2$, and for $m$ odd the $\Ga$-factor has the form:
 \\ $   
\tilde{\Gamma}(k,s)  =  i^{mk} 2^{-m(k+1)} \pi^{-m(s+k)}                           
                            {\Gamma}_m(k+s)$.
{}
\end{exam}
 {\sf {} We use this series with} $\bpsi=\chi^2\psi_1\psi_2\overline\psi_3$, 
$k=k_2+k_3-k_1\ge 2$, 
$m=3$, $\ds\kappa=\frac {m+1} 2=2$, $[m/2]=1$.



\begin{theo}[Siegel, Shimura \cite{Shi83}, P. Feit \cite{Fei86}]
\label{sd-5.1}
 Let $m$ 
be an {\sf odd integer} such that $2k>m$, and $N>1$ be an integer, then:

  For an integer $s$  such that 
{
$s=-r $, $0\leq r \leq k-\kp$}, {}
  there is the
  following Fourier expansion {}
  {\begin{equation}  \label{EQposcoefofGplus}
      G^{\star}({z}, -r) =  G^{\star}({z}, -r; k, {\bpsi}, N)
= \sum_{A_m \ni {\scrT} \ge 0} 
a({\scrT}, R)
q^{\scrT} ,
    \end{equation}}{}
  where for $s > (m+2-2k)/4$ in (\ref{EQposcoefofGplus}) 
    the only non-zero terms
   occur for positive definite ${\scrT}>0$,
\end{theo}
\subsubsection*{Fourier coefficients of Siegel-Eisenstein series (continued)}

{\it

  \begin{align}\label{aTR}&
a({\scrT}, R)=
  M({\scrT},{\bpsi},k-2r) \cdot 
\det({\scrT})^{k-2r-\kp}
Q(R,{\scrT}; k-2r, r), 
\\ & \label{FinEu}
M ({\scrT},k-2r,{\bpsi})=\prod\limits_{\ell\vert\det(2{\scrT})}
M _\ell({\scrT}, {\bpsi}(\ell)\ell^{-k+2r})
\end{align}
{}
polynomials
$
Q(R,{\scrT}; k-2r, r)
$ are given by  (\ref{QRT}), and for all
  ${\scrT}>0$, ${\scrT} \in A_m$, 
 is a finite Euler
product, in which
$M _\ell({\scrT}, x)\in\dbZ[x]$. 
 \qed
}


\section{Statement of the Main Result}
\subsubsection*{Main Theorem (on $p$-adic analytic function 
 in four variables)}

{\it
{\sf 1)} The function
$\ds
\Lr_{\textbf f}: (s, k{}_1, k{}_2, k{}_3)\mapsto \frac {\ang{\textbf f^{\,0}, {\cal E}(- r, \chi)},}
{ \ang{\textbf f^{\,0}, \textbf f_0},}
$
depends
{ $p$-adic analytically
on four variables} $(\chi\cdot y_p^r, k{}_1, k{}_2, k{}_3)
\in X\times \Br_1\times\Br_2\times\Br_3$;

{\sf 2)}
 { Comparison  of complex and $p$-adic values:} for all $({k_1, k_2, k_3})$ 
in an affinoid neighborhood 
$
\Br= \Br_1\times\Br_2\times\Br_3\subset X^3, 
$
 satisfying $k{}_1 \le k_2{} + k_3{} -2$:
 the values
at $s=k_2 + k_3 -2-r$ 
coincide with the normalized critical special values
\begin{align}\label{L^*}&
{}
L^*({f_{1,k_1{}}}\otimes {f_{2,k_2{}}}
\otimes {f_{3, k_3{}}}, k_2{}+k_3{}-2-r, \chi)  \\
\nonumber & (r=0, \cdots, k_2{}+k_3{}-k{}_1-2), 
\end{align}
 for  Dirichlet characters
$\chi\bmod Np^v, v\ge 1$,}
{}
such that all three corresponding Dirichlet characters $\chi_j$ 
have  $Np$-complete conductors:

\subsubsection*{Main Theorem (continued)} 

{\it 
\begin{align}\label{EQchij}&
\chi_1\bmod Np^v=\chi, \ 
\chi_2\bmod Np^v=\psi_2\bar\psi_3\chi, \\ &\nonumber 
\chi_3\bmod Np^v=\psi_1\bar\psi_3\chi, 
\bpsi=\chi^2\psi_1\psi_2\overline\psi_3. 
\end{align}
{}
The normalisation of $L^*$ in (\ref{L^*}) is the same as in Theorem C below.

{\sf 3)} {{} Dependence on $x\in X$:} let $H= [2{\rm ord}_p(\lambda)]+1$.
For any fixed $({k{}_1, k{}_2, k{}_3})\in\Br$ and $x=\chi\cdot y_p^r$
the function 
$$
x\longmapsto 
\frac {\ang{\textbf f^{\,0}, {\cal E}(- r, \chi)},}{ \ang{\textbf f^{\,0}, \textbf f_0},}
$$
extends to a $p$-adic analytic function of type $o(\log^H(\cdot))$ of the variable
 $x\in X$.
}

\rem.
The function
$\Lr_{\textbf f}$ depends on the variables $(s, k_1, k_2, k_3)$ in a different way:
it is a mixture of the $p$-adic  Mellin transform (in $s$), and of a rigid analytic function 
(in $k_1, k_2, k_3)$. 


\subsubsection*{Outline of the proof} 

{\sf 1)} 
$\bullet$
At each classical weight $(k{}_1, k{}_2, k{}_3)$ let us use the equality {{}
$$
{\ang{\textbf f^{\,0}, {\cal E}(- r, \chi)},}={\ang{\textbf f^{\,0}, 
\pi_\lambda({\cal E}(- r, \chi)}),}
$$
} 
which is
deduced from
the definition of the projector $\pi_{\lambda}$:
$
\Ker \pi_{\lambda}:=\bigcap_{n\ge 1}\Im 
(U_T-{\lambda} I)^n, 
\Im \pi_{\lambda}:=\bigcup_{n\ge 1}\Ker 
(U_T-{\lambda} I)^n. 
$

{}{{}
Notice that the coefficients of 
 {{} $
{\cal E}(- r, \chi)
 \in \scrM (\scrA)  
$
}
 depend $p$-adic analytically on $({k_1, k_2, k_3})\in\Br=\Br_1\times\Br_2\times\Br_3$}, 
where ${\scr A}= {\scr A}(\Br_1\times\Br_2\times\Br_3)$ is the $p$-adic Banach algebra
of rigid-analytic functions on  $\Br$.

\subsubsection*{ Interpolation to all $p$-adic weights:}

{$\bullet$} At each classical weight $(k{}_1, k{}_2, k{}_3)$ the scalar product {{}
$\ang{\textbf f^{\,0}, {\cal E}(- r, \chi)},$} is given by the first coordinate 
of $\pi_\lambda(
{\cal E}(- r, \chi))$ with respect to an orthogonal basis of $\scrM^\lambda(\scrA)$
containing $\textbf f_0$ with respect to {\sf{} Hida's algebraic Petersson product 
$\ang{g, h},_a:=\ang{g^\rho|\mat 0, -1, Np, 0, , h},$}, see \cite{Hi90}.


Let us extend the linear form
$\ell(h)=
\frac {\ang{\textbf f^{\,0}, h},}{ \ang{\textbf f^{\,0}, \textbf f_0},}$
(defined for classical weights),
 to  Coleman's type submodule 
of overconvergent families $h\in \scrM^\lambda(\scrA)^\dagger\subset \scrM^\lambda(\scrA)$ as
{{}
the first coordinate of $h$ with respect to some
$\Ar$-basis of 
eigenfunctions of all (triple) Hecke operators $T_q$ for $q\nmid Np$, 
having the first basis vector ${\textbf f}_0\in \scrM^\lambda(\scrA)^\dagger$.} 

The linear form $\ell$ can be characterized as a normalized eigenfunction 
of the adjoint Atkin's operator, acting on the dual $\Ar$-module of $\scrM^\lambda(\scrA)^\dagger$:
$\ell({\textbf f}_0)=1$.

In order to extend $\ell$ to 
$h={\cal E}(- r, \chi)$, 
we need to choose a certain representative  of ${\cal E}(- r, \chi)$ in 
the $\Ar$-submodule $\scrM^\lambda(\scrA)^\dagger$ , which is locally free of finite rank.

\subsubsection*{A representative  of ${\cal E}(- r, \chi)$ in 
the  (locally free of finite rank $\Ar$-submodule) $\scrM^\lambda(\scrA)^\dagger$} 

Choose a (local) basis $\ell^1, \cdots, \ell^n$ given by some {\it{} triple Fourier coefficients}
of the dual (locally free of finite rank) $\Ar$-module $\scrM^\lambda(\scrA)^\dagger{}^*$.

Then
define
$$
\ell=\beta_1\ell^1+ \cdots+ \beta_n\ell^n,
$$
 where 
{{} $\beta_i=\ell(\ell_i)\in\Ar$}, and $\ell_i$ denotes the dual basis of  
$\scrM^\lambda(\scrA)^\dagger$: $\ell^j(\ell_i)=\delta_{ij}$. 
At each $p$-adic weight $(k{}_1, k{}_2, k{}_3)\in\Br$ let us define {}
$$
\ell({\cal E}(- r, \chi)):=\beta_1\ell^1({\cal E}(- r, \chi))+ \cdots+ \beta_n\ell^n({\cal E}(- r, \chi))
{{} \mbox{ (belongs to } \Ar)},
$$
{}
where $\beta_i=\ell(\ell_i)\in \Ar$, and $\ell^i({\cal E}(- r, \chi))\in \Ar$ 
are certain Fourier coefficients of the seies ${\cal E}(- r, \chi)$.

\subsubsection*{Conclusion}


There exists an element 
$$
\tilde{\cal E}(- r, \chi)\in \scrM^\lambda(\scrA)^\dagger\subset\scrM(\scrA)^\dagger
$$
such that 
\fbox{
$
\ell({\cal E}(- r, \chi))=\ell(\tilde{\cal E}(- r, \chi))
$
}
(at each weight $(k{}_1, k{}_2, k{}_3)$).
{} 
In fact, let us define
$$
{}
\tilde{\cal E}(- r, \chi):=\ell^1({\cal E}(- r, \chi))\ell_1+ \cdots+ \ell^n({\cal E}(- r, \chi))\ell_n
$$
$
\Rightarrow
\ell(\tilde{\cal E}(- r, \chi))=\ell(\ell_1)\ell^1({\cal E}(- r, \chi))+ \cdots+ \ell(\ell_n)\ell^n({\cal E}(- r, \chi))
$\\
$
=\beta_1\ell^1({\cal E}(- r, \chi))+ \cdots+ \beta_n\ell^n({\cal E}(- r, \chi))
$\\
$
=\ell({\cal E}(- r, \chi))
$ (at each 
weight $(k{}_1, k{}_2, k{}_3)$).

Thus,  the dependence of $\ell({\cal E}(- r, \chi))\in \Ar$
on $(k_1, k_2, k_3)\in X^3$ is $p$-adic analytic.

\

In order  to prove the remaining statements {\sf 2)},  {\sf 3)},   the dependence on $x=\chi\cdot y_p^r$ 
is studied
in the next section.











\section{Distributions and admissible measures}
\label{sec:adm}
\subsubsection*{Distributions and measures with values in Banach modules}

\subsection*{Notation}
\begin{tabular}{ll}
 \begin{tabular}{l}
$\Ar$ \\
$V$\\ 
${\cal C}(Y,\Ar)$\\ 
$\cup$ \\
${\cal C}^{loc-const}(Y,\Ar)$\\ 
\end{tabular} 
 \begin{tabular}{l}\\
(a $p$-adic Banach algebra) \\
(an $\Ar$-module)\\ 
(the $\Ar$-Banach algebra \\ of {\sf{} continuous functions}
 on $Y$ ) 
\\ 
(the $\Ar$-algebra \\ of {\sf{} locally constant functions}
 on $Y$ ) 

\end{tabular}
\end{tabular}

\begin{defi}[Distributions and measures]\label{defi:distr}

{a)}
A {\sf{} distribution} $\Dr$ on $Y$ with values in $V$ is 
an $\Ar$-linear form 
$$
\Dr: {\cal C}^{loc-const}(Y,\Ar)\to V, \ \ \varphi\mapsto\Dr(\varphi)=\int_Y \varphi d \Dr.
$$

{b)}
A {\sf{} measure} $\mu$ on $Y$ with values in $V$ is 
a continuous  $\Ar$-linear form 
$$
\mu: {\cal C}(Y,\Ar)\to V, \ \ \varphi\mapsto\int_Y \varphi d \mu.
$$

\end{defi}

\

The integral $\ds \int_Y \varphi d\mu$ can be defined for 
any continuous function $ \varphi$, and any
 bounded distribution $\mu$, 
using the Riemann sums.


\subsection*{Admissible measures of Amice-Vélu}

\subsubsection*{Admissible measures}

Let $h$ be a positive integer.
A more delicate notion of an $h$-admissible measure
 was introduced
 in \cite{Am-V} by Y. Amice, J.  Vélu  (see also \cite{MTT}, 
\cite{V}):

\begin{defi}\label{defi:adm}\ 

{a)}
For $h\in \N, h\ge 1$ let ${\scrP}^h(Y, \Ar)$ denote
 the $\Ar$-module of {\sf{} locally polynomial functions} of degree $<h$
 of the variable $y_p:Y\to \Z_p^\times \hookrightarrow \Ar^\times$;
 in particular, 
$$
{\scrP}^1(Y, \Ar)  = {\cal C}^{loc-const}(Y, \Ar)
$$
 (the 
 $\Ar$-submodule of {\sf{} locally constant functions}).
Let also denote
 ${\cal C}^{loc-an}(Y, \Ar)$ the 
 $\Ar$-module of {\sf{} locally analytic functions}, so that
  $$
{\scrP}^{1}(Y, \Ar)\subset  {\scrP}^h(Y, \Ar) \subset 
 {\cal C}^{loc-an}(Y, \Ar)\subset  {\cal C}(Y, \Ar).
$$
\end{defi}




{\it
{b)}
Let $V$ be a normed $\Ar$-module with the norm $|\cdot|_{p,V}$. 
For a given positive integer $h$ an $h$-{\it admissible measure} 
on $Y$ with values in $V$ is an $\Ar$-module homomorphism 
$$
\tilde \Phi :  {\scrP}^h(Y, \Ar)\to V
$$
such that for fixed $a\in Y$ and for ${v}\to \infty$
the following {\sf {} growth condition} is satisfied: 
\begin{eqnarray}\label{equation:adm} & & 
\left |\int_{a+(Np^{v})}(y_p-a_p)^{h^\prime}d\tilde \Phi\right |_{p,V}
= o(p^{-{v}(h^\prime-h)})\\ & &  {\rm for \ \ all}
\ \ \ h^\prime =0, 1, \dots , h-1, a_p:=y_p(a) \nonumber
\end{eqnarray} 
}

The condition (\ref{equation:adm}) allows 
 to  integrate {\sf{} the locally-analytic functions
on $Y$ along $\tilde \Phi$
using Taylor's expansions!
}
This means: there exists a unique extension of $\tilde \Phi$ to 
 ${\cal C}^{loc-an}(Y, \Ar)\to V$. 





\subsection{$U_p$--Operator and the method of canonical projection}

\subsubsection*{Using the canonical projection $\pi_{\lambda}$}

 We construct
our $H$-admissible measure 
$\widetilde \Phi^{\lambda}: {{\scrP}}^H(Y, {\Ar})\to{\scrM}({\Ar})$ 
out of a sequence of distributions
$
\Phi_{r} : {{\scrP}}^1(Y, {\Ar})\to {\scrM} ({\Ar})
$ defined {\sf{} on local monomials $y_p^r$  of each degree $r$} by the rule
$$
{}
\int_{Y } \chi y_{p}^r \,d\widetilde \Phi^{\lambda}
=\pi_{\lambda}(\tilde{\cal E}(- r, \chi)), 
\mbox{ where } \tilde{\cal E}(- r, \chi)\in M={\scrM}({\Ar}).
$$
Here 
$\tilde{\cal E}(- r, \chi)$ takes values in an  ${\Ar}$-module 
{}
$$
M={\scrM}({\Ar})\subset {\Ar}[\![q_1, q_2, q_3]\!][R_1, R_2, R_3]
$$
{}
 of nearly holomorphic (overconvergent) triple modular
forms over ${\Ar}$ 
(for $0\le r\le H-1$, $H=[2\ord_p \lambda_p]+1)$, 
and the formal series $\tilde{\cal E}(- r, \chi)$
 was constructed in the proof of 1) of Main Theorem.

\subsubsection*{Definition of the canonical projection $\pi_{\lambda}$}

Here
$\Ar$ is an ${\mathbb C}_p$-algebra, and
{}
$
{\lambda}\in {\Ar}^\times 
$
{}
 is a fixed non-zero eigenvalue 
of triple Atkin's operator
 $U_T= U_{T, p}$, acting on ${\scrM}({\Ar})$, 
 $$\pi_{\lambda}: {\scrM}({\Ar})\to {\scrM}({\Ar})^{\lambda}
$$
{}
is the {\sf{} canonical projection
 operator}
onto the maximal ${\Ar}$-submodule 
${\scrM}({\Ar})^{\lambda}$
over which the operator $U_T-{\lambda} I$ is nilpotent
(we call ${\scrM}({\Ar})^{\lambda}$ the $\lambda$-characteristic 
submodule of ${\scrM}({\Ar})$).\\
The projector $\pi_{\lambda}$ is defined by its kernel:
{
$
\Ker \pi_{\lambda}:=\bigcap_{n\ge 1}\Im 
(U_T-{\lambda} I)^n, 
$
$
\Im \pi_{\lambda}:=\bigcup_{n\ge 1}\Ker
(U_T-{\lambda} I)^n.
$
}


\section{Triple modular forms} 

Triple modular forms
are certain formal series
{
\begin{align*}&
g=\sum\limits_{n_1, n_2, n_3=0}^\infty a(n_1, n_2, n_3; R_1, R_2, R_3)q_1^{n_1}q_2^{n_2}q_3^{n_3}
\\ & \in 
{\Ar}[\![q_1, q_2, q_3]\!][R_1, R_2, R_3], \mbox{ where }z_j=x_j+iy_j\in \dbH, \ R_j=(4\pi y_j)^{-1}, 
\end{align*}
}
with the property that for ${\Ar}={\mathbb C}$, 
  the series converges to a $\scrC^\infty$-modular form on $\dbH^3$ of a given
 weight $(k_1,  k_2, k_3)$ and character $(\psi_1, \psi_2, \psi_3)$, $j=1, 2, 3$.

The coefficients $a(n_1, n_2, n_3; R_1, R_2, R_3)$
 are polynomials in ${\Ar}[R_1, R_2, R_3]$, and
the triple Atkin's operator 
 is given by
\\
{
$
U_T(g)=\sum\limits_{n_1, n_2, n_3=0}^\infty a(pn_1, pn_2, pn_3; 
pR_1, pR_2, pR_3)q_1^{n_1}q_2^{n_2}q_3^{n_3}.
$
}


\subsection*{Eigenfunctions of $U_p$ and of $U_p^*$.}
\subsubsection*{ Functions $f_{j, 0}$ and $f_{j}^{\,0}$}
Recall that  
for any primitive
 cusp eigenform 
$f_j = \sum_{n=1}^\infty a_n(f)q^n$,    
there is an eigenfunction  
{{}
$f_{j, 0} = \sum_{n=1}^\infty a_n(f_{j, 0})q^n \in \Qb[\![q]\!]$
}{}
 of $U=U_p$
with the eigenvalue
 $\alpha=\alpha_{ p,j}^{(1)}\in \Qb$
($U(f_0)=\alpha f_0$) given by
{{}
\begin{align}\label{EQf_0}&
f_{j, 0}=f_j-\alpha_{ p,j}^{(2)}f_j|V_p=f_j-
\alpha_{ p,j}^{(2)}p^{-k/2}f_j|\mat p, 0, 0, 1, 
\\ &\nonumber
\sum_{n=1}^\infty a_n(f_{j, 0})n^{-s}= 
\sum_{\scriptstyle n=1\atop \scriptstyle {p\nmid n}}^\infty
  a_n(f_j)n^{-s}(1-\alpha_{ p,j}^{(1)}p^{-s})^{-1}.
\end{align}
}


Moreover, there is an eigenfunction 
{
$ f_j^{\,0}$ of $U_p^*$
}
 given by
\begin{align}\label{2tv-(0.4)}
  f_j^{\,0} = 
 f ^{\rho}_{j, 0} \Big \vert _k  \ 
\mat 0,-1,  Np, 0, , \mbox{ where }
f^{\rho}_{j, 0}=\sum_{n=1}^\infty 
\overline {a(n, f_0 )} q^n.
\end{align}
Therefore, 
$U_T(f_{1, 0}\otimes f_{2, 0}\otimes f_{3, 0})=\lambda(f_{1, 0}\otimes f_{2, 0}\otimes f_{3, 0})$.
{}


\section{Critical values
of the $L$ function $L(f_1\otimes f_2\otimes f_3, s, \chi )$ 
} 

\subsubsection*{Choice of Dirichlet characters} 


For  an arbitrary Dirichlet character $\chi\bmod Np^v$
 consider the following Dirichlet characters:
\begin{align}  & 
\chi_1\bmod Np^v=\chi, \ 
\chi_2\bmod Np^v=\psi_2\bar\psi_3\chi, \\ &\nonumber 
\chi_3\bmod Np^v=\psi_1\bar\psi_3\chi, 
\bpsi=\chi^2\psi_1\psi_2\overline\psi_3; 
\end{align}
later on we
impose the condition that 
the conductors of 
the corresponding primitive characters
$
\chi _{0, 1}, \chi_{0, 2}, \chi_{0, 3}
$
are $Np$-completes (i.e. have the same prime divisors as resp. those of $Np$).

\subsubsection*{{\sc {Theorem A}} 
{\sc (algebraic properties of the triple product)}
}
{\it
Assume that $k_2+k_3-k_1\ge 2$, then
for all pairs $(\chi, r)$ such that 
the corresonding Dirichlet characters $\chi_j$ 
 have  $Np$-complete conductors, and 
  $0 \leq  r \leq k_2+k_3- k_1- 2$,  we have that
{}
$$ 
\frac{\La( f_1^{\rho}\otimes f_2^{\rho}\otimes f_3^{\rho}, k_2+k_3-2-r,
\psi_1\psi_2\chi )}
{\langle f_1^{\rho}\otimes f_2^{\rho}\otimes f_3^{\rho},  
f_1^{\rho}\otimes f_2^{\rho}\otimes f_3^{\rho}\rangle_T} 
 \in \Qb 
$$
{}
  where 
\begin{align*}&
{\langle f_1^{\rho}\otimes f_2^{\rho}\otimes f_3^{\rho},  
f_1^{\rho}\otimes f_2^{\rho}\otimes f_3^{\rho}\rangle_T}
:=\langle f_1^{\rho}, f_1^{\rho}\rangle_N\langle f_2^{\rho}, f_2^{\rho}\rangle_N \langle f_3^{\rho},  
f_3^{\rho}\rangle_N\\ &\qquad
=\langle f_1, f_1\rangle_N\langle f_2, f_2\rangle_N \langle f_3,  
f_3\rangle_N.
\end{align*}

}

\section{Theorems B-D} 
\subsubsection*{Recall: the $p$-adic weight space and the Mellin transform}
The {\sf$p$-adic weight space} is
the group
  $X={\rm Hom}_{cont}(Y,\cb_p^\times)$ 
of (continuouos) $p$-adic characters 
of  the commutative profinite group  
$\displaystyle{ Y 
= \lim_{\buildrel \longleftarrow\over v } (\zb/Np^v\zb)^*}$ 

The group  $X$ is isomorphic to a finite union of discs 
$U=\{z\in\cb_p\ |\ |z|_p<1\}$.

A $p$-adic $L$-function  $L_{(p)}:X\to \cb_p$ is 
a certain meromorphic function on
 $X$.  
Such a function usually come from a $p$-adic measure $\mu$ on  $Y$
({\it bounded} or {\it admissible} in the sense of Amice-Vélu, see
\cite{Am-V}).
The {\sf{} $p$-adic Mellin transform} of $\mu$
is given 
for all $x\in X$ by
$$
L_{(p)}(x)=\int_{Y_{N,p}}x(y)\dd\mu(y), L_{(p)}: X\to {\mathbb C}_p
$$

\subsubsection*{
 Theorem B (on admissible measures attached to
 the triple product:fixed balanced weights case)}
{\it
  Under the assumptions as above 
there
  exist
a ${\mathbb C}_p$-valued measure
$\tilde\mu^\lambda_{f_1\otimes  f_2\otimes   f_3}$ on $Y_{N,p}$,
 and  a ${\mathbb C}_p$-analytic function
\\ $
 \Dr_{(p)}(x, f_1\otimes f_2\otimes f_3) : X_p \rightarrow {\mathbb C}_p,  
$
given for all $x\in X_{N,p}$ by the integral 
 $
\Dr_{(p)}(x, f_1\otimes f_2\otimes f_3)
=\int_{Y_{N,p}}x(y)\dd\tilde\mu^\lambda_{f_1\otimes  f_2\otimes   f_3}(y), 
$
and having the following properties:
\\
}

{\it
{\rm (i)} for all pairs $(r, \chi)$ such that 
$\chi \in X_{N, p}^{\tors}$, 
and
all 
three corresonding Dirichlet characters $\chi_j$ 
have $Np$-complete conductor 
$(j=1, 2, 3)$, and ${r}\in \ZZ$ 
is an integer with
  $0 \leq  r \leq k_2+k_3-k_1 - 2$,  
the following equality holds:
\begin{align*}   
\Dr_{(p)}&(\chi x_p^{r}, f_1\otimes f_2\otimes f_3) =
 i_p \Big( 
\frac{
(\psi_1
\psi_2)(2)
C_{\chi}^
       {4(k_2+k_3-2-r)}}
{G(\chi_1)G(\chi_2)G(\chi_3) 
G(\psi_1
\psi_2\chi_1)\lambda_p^{2v}} 
\end{align*}
$\ds
\frac{\La (f_1^{\rho}\otimes f_2^{\rho}\otimes f_3^{\rho}, k_2+k_3-2-{r}, \psi_1\psi_2\chi)}
       {\langle f_{1}^{0}\otimes f_{2}^{0}\otimes f_{3}^{0},  
f_{1, 0}\otimes f_{2, 0}\otimes f_{3, 0}\rangle_{T, Np} }
\Big)  
$
\\ 
where $v=\ord_p(C_{\chi})$, 
$G(\chi)$ denotes  the Gau\ss\ sum 
 of a primitive 
Dirichlet character $\chi_0$ attached to $\chi$
 (modulo the conductor of 
 $\chi_0$), 
}

{\it 
{\rm (ii)} if $\ord_{p}\lambda_p = 0$ then the holomorphic function in (i)
 is a {\sf{} bounded ${\mathbb C}_{p}$-analytic function};
\\

{\rm (iii)} in the general case (but assuming that
 $\lambda_p\neq 0$)
  the holomorphic function in (i) belongs to the {\sf{} type
  $o(\log(x^{H}_{p}))$ with $H = [2 \ord_p \lambda_p] + 1$} 
and it can be
  represented as the Mellin transform  of 
the $H$-admissible ${\mathbb C}_p$-valued measure
$\tilde\mu^\lambda_{f_1\otimes  f_2\otimes   f_3}$ 
(in the sense of Amice-Vélu) on $Y$ 
\\ 
{\rm (iv)} Let $k=k_2+k_3-k_1\ge 2$. 
If $H \leq k - 2$ then the function  $\Dr_{(p)}$
 is uniquely
  determined by the above conditions (i).
}


Let us  describe now 
$p$-adic measures
attached to Garrett's triple product of three Coleman's families
\begin{align}\label{EQfj'}
k_j{}\mapsto \{f_{j,k_j{}}= \sum_{n=1}^\infty a_{n,j}(k{})q^n\} (j=1, 2, 3).
\end{align} 
Consider the product of three eigenvalues:
$$
\lambda= \lambda_p(k_1{}, k_2{}, k_3{})
= \alpha_{p, 1}^{(1)}(k{}_1)\alpha_{p, 2}^{(1)}(k{}_2)\alpha_{p, 3}^{(1)}(k_3{})
$$
and assume that the slope of this product  
$$
\sigma= \ord_p(\lambda(k_1{}, k_2{}, k_3{}))=
\sigma(k_1{}, k_2{}, k_3{})=
\sigma_1+\sigma_2+\sigma_3
$$
 is {\sf{} constant and positive} for all triplets $(k_1{}, k_2{}, k_3{})$ 
 in an appropriate $p$-adic neighbourhood
 of the fixed triplet of weights $(k_1, k_2, k_3)$.


Let $\Ar=\Ar({\cal  B})$ denote  an affinoid algebra $\Ar$
associated
 with an analytic space 
${\cal  B}={\cal  B}_1\times{\cal  B}_2\times{\cal  B}_3$, 
a  neighbourhood
 around $(k_1, k_2, k_3)\in {{X}}^3$ (with
 a given  $k$ and  $\psi \bmod N$).


\subsubsection*{Theorem C (on $p$-adic measures for 
families of triple products) }{\it\ 
Put $H= [2{\rm ord}_p(\lambda)]+1$.
There exists a sequence of distributions 
$\Phi_r$
on $Y$  with values in  $\mc =\mc(\Ar)$ 
giving 
 an $H$-admissible measure $\tilde\Phi^\lambda$
 with values in
  $\mc^\lambda\subset \mc $
with the following properties:

There exists an $\Ar$-linear form $\ell = 
\ell_{{\textbf f}_1\otimes{\textbf f}_2\otimes{\textbf f}_3, 
\lambda}: \Mc(\Ar)^\lambda\to \Ar
$ (given by   (\ref{ell}), 
such that the composition 
$$
\tilde\mu
= \tilde\mu_{{\textbf f}_1\otimes{\textbf f}_2\otimes{\textbf f}_3, \lambda} 
:=\ell_{{\textbf f}_1\otimes{\textbf f}_2\otimes{\textbf f}_3, \lambda}(\tilde
\Phi^{\lambda} )
$$
is an $H$-admissible measure
with values in $\Ar$, 
and
for all $({k{}_1, k{}_2, k{}_3})$ 
in the affinoid neighborhood $
\Br= \Br_1\times\Br_2\times\Br_3, 
$ as above, 
 satisfying $k{}_1 \le k_2{} + k_3{} -2$
}

{\it 
we have that
 the evaluated integrals
$$
ev_{(k{}_1, k{}_2, k{}_3)}\left((\ell_{ {\textbf f}_1\otimes{\textbf f}_2\otimes{\textbf f}_3,
 \lambda})
(\tilde\Phi^\lambda)(y_p^r\chi )\right)
$$ 
on the  arithmetical chracters $ y_p^r\chi$
coincide with the critical special values
{}
$$
\La^*({f_{1,k_1{}}}\otimes {f_{2,k_2{}}}
\otimes {f_{3, k_3{}}}, k_2{}+k_3{}-2-r, \chi)  \ \ 
$$
{}
for $r=0, \cdots, k_2{}+k_3{}-k{}_1-2$, 
 and for all Dirichlet characters
$\chi\bmod Np^v, v\ge 1$,
with all three corresonding Dirichlet characters $\chi_j$ 
given by (\ref{EQchij}),
having  $Np$-complete conductors.
Here
the normalisation of $\La^*$ includes at the same time
certain Gauss sums, Petersson scalar products, 
powers of $\pi$ and of $\lambda(k{}_1, k{}_2, k{}_3)$, 
and a certain finite Euler product.
}

\subsubsection*{The $p$-adic Mellin transform and four
variable $p$-adic analytic functions}

Any $h$-admissible measure $\tilde\mu$ on $Y$ with values in 
a $p$-adic Banach algebra $\Ar$ 
 can be caracterized
its {\it Mellin transform} 
 $\Lr_{\tilde\mu} (x)$ 
$
\Lr_{\tilde\mu}:X\to \Ar$, defined by 
$
\Lr_{\tilde\mu} (x) = \int_{Y} x (y)d\tilde\mu (y), 
$ 
where 
$  x \in X, \  \Lr_{\tilde\mu} (x)\in \Ar, 
$

\noindent
{\sf  Key property of $h$-admissible measures $\tilde \mu$:}
 its Mellin transform $\Lr_{\tilde\mu}$ 
is {\sf{} analytic} with values in $\Ar$.




Let 
$
\Ar=\Ar({\cal  B})=\Ar_1\hat\otimes\Ar_2\hat\otimes\Ar_3
=
\Ar({\cal  B}_1)\hat\otimes\Ar({\cal  B}_2)\hat\otimes\Ar({\cal  B}_3)
$
{}
 denote again the Banach algebra $\Ar$
where ${\cal  B}$ is an affinoid neighbourhood
 around $(k_1, k_2, k_3)\in {{X}}^3$ (with
 a given  integer $k$ and Dirichlet character $\psi \bmod N$).
\subsubsection*{
Theorem  D (on $p$-adic analytic function 
 in four variables)}

{\it\ 
Put $H= [2{\rm ord}_p(\lambda)]+1$.
There exists a
{ $p$-adic analytic function 
 in four variables} $(x, {\textbf s})
\in X\times \Br_1\times\Br_2\times\Br_3\subset X\times X\times X\times X$:
$$
\Lr_{\tilde\mu}: (x, {\textbf s}) \longmapsto 
ev_{\textbf s}(\Lr_{\tilde\mu(x)})  \ \ \
 (x \in X, \ \ \Lr_{\tilde\mu} (x)\in \Ar).
$$
with values in ${\mathbb C}_p$, 
such that 
for all $({k{}_1, k{}_2, k{}_3})$ 
in the affinoid neighborhood as above 
$
\Br= \Br_1\times\Br_2\times\Br_3, 
$
 satisfying $k{}_1 \le k_2{} + k_3{} -2$, 
we have that
 the special values
$
\Lr_{\tilde\mu} (x, {\textbf s}) 
$
at the  arithmetical chracters $ x=y_p^r\chi$, 
and ${\textbf s}= ({k{}_1, k{}_2, k{}_3})\in\Br$
coincide with the normalized critical special values
$$
L^*({f_{1,k_1{}}}\otimes {f_{2,k_2{}}}
\otimes {f_{3, k_3{}}}, k_2{}+k_3{}-2-r, \chi)  \ \ (r=0, \cdots, k_2{}+k_3{}-k{}_1-2), 
$$
 for  Dirichlet characters
$\chi\bmod Np^v, v\ge 1$,
such that all three corresonding Dirichlet characters $\chi_j$ 
given by (\ref{EQchij}),
have  $Np$-complete conductors
where the same normalisation of $L^*$ as in Theorem C.
}

{\it
Moreover, for any fixed ${\textbf s}= ({k{}_1, k{}_2, k{}_3})\in\Br$
the function 
$$
x\longmapsto \Lr_{\tilde\mu}(x, {\textbf s}) 
$$
is $p$-adic analytic of type $o(\log^H(\cdot))$.
}

\

Indeed,  we obtain the  function in question
$\Lr_\mu (x, {\textbf s})$
by evaluation at 
$$
{\textbf s}=((s_1, \psi_1),  (s_2, \psi_2),  (s_3, \psi_3))\in \Br:
$$
this is  a {\sf{} $p$-adic analytic function 
 in four variables} $(x, {\textbf s})
\in X\times \Br_1\times\Br_2\times\Br_3\subset X\times X\times X\times X$:
$${{}
\Lr_{\tilde\mu} (x, {\textbf s}) := 
ev_{\textbf s}(\Lr_{\tilde\mu}) (x) \ \ \
 (x \in X, \ 
{\textbf s}
\in  \Br_1\times\Br_2\times\Br_3, 
\  \Lr_{\tilde\mu} (x)\in \Ar).}
$$



{\sf  This is  a joint work in progress with S.Boecherer},  
 we use: 

1) the {\sf{} higher twists of the Siegel-Eisenstein series}, 
introduced in \cite{B\"o-Sch}, 

2)
  {\sf{} Ibukiyama's differential operators} (see \cite{Ibu}, 
\cite{BSY}).


\section{Ideas of the Proof}\label{SchPr}
\subsection
{Boecherer's higher twist} 
We define the higher twist of the series
$F_{\chi,r}=\sum_{{\cal T}}
a_{\chi,r}(R, {\cal T})q^{\cal T}
$ by some Dirichlet  characters $\bar \chi_1, \bar\chi_2, \bar\chi_3$ 
as the following {{\sf{} formal nearly holomorphic Fourier expansion}}: 
\begin{align} \label{Fchichi} 
F_{\chi,r}= 
\sum_{{\cal T}}
\bar\chi_1(t_{12})
\bar\chi_2(t_{13}) 
\bar\chi_3(t_{23})
a_{\chi,r}(R, {\cal T})q^{\cal T}.
\end{align}
The seies (\ref{Fchichi}) is a Siegel modular form of some higher level, 
but it has additional symmetries with respect to symplectic embedding
$\iota_3: \Gamma_0(Np^{2v})\times\Gamma_0(Np^{2v})\times\Gamma_0(Np^{2v})\to\Sp_3$:
{{} its triple Nebentypus character does not depend on $\chi\bmod Np^v$, 
and is equal to $(\psi_1, \psi_2, \psi_3)$,  
}
if we choose Dirichlet characters as in (\ref{EQchij}):
\begin{align*}&
\chi_1\bmod Np^v=\chi, \ 
\chi_2\bmod Np^v=\psi_2\bar\psi_3\chi, \\ &\nonumber 
\chi_3\bmod Np^v=\psi_1\bar\psi_3\chi, 
\bpsi=\chi^2\psi_1\psi_2\overline\psi_3.
\end{align*}

We use the  Siegel-Eisenstein series 
$
F_{\chi, r}
$
which depends on the character $\chi$, 
but its precise nebentypus character is
$\bpsi=\chi^2\psi_1\psi_2\overline\psi_3$, 
and it is defined by 
$F_{\chi, r}=
G^{\star}({z}, -r; k, {(Np^v)}^2,\bpsi)$,
where 
${z}$ denotes a variable 
in the Siegel upper half space $\HH_3$, 
and the normalized series $G^{\star}({z}, s; k, {(Np^v)}^2,\bpsi)$ is given by (\ref{eqG*}).

This series depends on $s=-r$, and for the critical  values
at integral points
 $s\in\Z$ such that $2-k\le s\le 0$, 
 it  represents a ({\it nearly})  {\it holomorphic} Siegel modular form
in the sense of Shimura \cite{ShiAr}:
$$
F_{\chi, r}
=\sum_{{\cal T}}
\det({\cal T})^{k-2r-\kp}
Q(R,{\cal T}; k-2r, r)
a_{\chi,r}({\cal T}) q^{\cal T}. 
$$



\subsection
{Ibukiyama's differential operator} 
We use an algebraic version of
 {\sf Ibukiyama's differential operator}, 
 which generalizes
the algebraic ``pull-back'':
 it transforms
 a nearly holomorphic {\sf Siegel modular form }
of weight $k{}$ 
to a nearly holomorphic {\sf triple modular form }
of weight $(k_1{}, k_2{}, k_3{})$ $(k{}=k{}_2+ k_3{}-k_1{})$.

On a holomorphic Siegel  modular form 
$
F=\sum_{\cal T} a({\cal T}) q^{\cal T}, 
$
this operator has the form
\begin{align}\label{IbuOp}
{\cal L}^{\lambda{},\nu{}}_{k{}}(F)
=
\sum_{{\cal T}}    
 \scrP(k_1{}, k_2{}, k_3{}, 0, \scrT)
a({\cal T}) q_1^{t_{11}}q_2^{t_{22}}q_3^{t_{33}},  
\end{align}
where $\lambda{}=k_1-k_3\ge \mu=k_1-k_2\ge 0$, and 
$\scrP(k_1{}, k_2{}, k_3{}; r; \scrT)$
 is certain Ibukiyama's polynomial,
equal to $(t_{11}t_{22}t_{33})^{\lambda{}}
(t_{12}t_{13}t_{23})^{\mu{}}$, if $r=0$.


As a result we obtain
a sequence of triple modular distributions $\Phi_r(\chi)$
 with values in the  tensor product 
{}
$
\Mr_T(\Ar)=
\Mr_{ }(\Ar)\widehat\otimes_\Ar\Mr_{}(\Ar)\widehat\otimes_\Ar
\Mr_{}(\Ar)
$
{}
 of three Banach $\Ar$-modules of arithmetical 
 nearly holomorphic modular forms
(the normalizing factor $2^r$ is neeeded in order to prove certain congruences
between $\Phi_{r}$). 
{\sf{} Note that 
${\Mr_T}(\Ar)$ is again a Banach $\Ar$-module
on which $U_T$ acts as a completely continuous operator.
}

{\sf The important property of the triple modular forms
$\Phi_{r}(\chi)$}: 
 the nebentypus character 
is fixed and is equal to $(\psi_1, \psi_2,  \psi_3)$
(for all $(k_1{}, k_2{}, k_3)$ and $\chi$ in question).

Using this property we compute the 
canonical projection $\pi_\lambda(\Phi_{r}(\chi))$
of the triple modular form $\Phi_{r}(\chi)$
   onto the $\lambda$-characteristic $\Ar$-submodule 
${\Mr_T}^{\lambda}(\Ar)
$
of the 
triple Atkin's operator $U_{T,p}$:
$$
\pi_\lambda:
{\Mr_T}(\Ar)\to \Mr_T^\lambda(\Ar).
$$

We prove 
that the resulting sequence of modular 
 distributions $\pi_\lambda(\Phi_{r})$ 
on the profinite group $Y$ 
 produces a certain $p$-adic admissible measure 
 $\tilde \Phi^\lambda$ (in the sense of Amice-Vélu,  \cite{Am-V})
 with values in 
 a certain {\sf{} locally free $\Ar$-submodule  of finite rank}
$$
\Mr_T^{\lambda}(\Ar)\subset {\Mr_T}(\Ar)
\subset
\mathop\bigcup_{v\ge 0} 
\mc_{T}(Np^v,\psi_1, \psi_2,\psi_3; \Ar)
$$
of formal nearly holomorphic triple modular forms of all
levels 
$Np^v$ and the fixed nebentypus characters
 $(\psi_1, \psi_2,\psi_3)$. 
We use
 congruences between triple 
 modular forms $\Phi_r(\chi)\in {\Mr_T}(\Ar)$
 (they have explicit  formal Fourier coefficients). 
 
Then we use  a {\sf {} general admissibility criterion}
 saying that these congruences 
imply $H$-admissibility
for their projections in 
$\Mr_T^{\lambda}(\Ar)$, 
where $H=[2\ord_p(\lambda)]+1$.

\subsection
{Algebraic linear form}
3)
{\sf From $\Mr_T^{\lambda}(\Ar)$ to $\Ar$}:
we use a $\Qb$-valued linear forms of type
$$
\scrL~: h\longmapsto \frac
{\ang{ f_1^{\,0}\otimes f_2^{\,0}\otimes  f_3^{\,0} , h},}
{\ang {f_1^{\,0},  f_{1,0}},  \ang {  f_2^{\,0},  f_{2,0}},  
\ang {  f_3^{\,0},  f_{3,0}}, }
$$
where $f_j^{\,0}$ is  the eigenfunction (\ref{EQf_0})  of 
the  conjugate  Atkin's operator $U_p^*$, 
and $f_{j,0}$  is the eigen\-function (\ref{2tv-(0.4)}) of
 $U_p$. 
The linear form $\scrL$ is
defined on the finite dimensional $\Qb$-vector characteristic subspace
\begin{align*}&
h\in \mc_{\textbf k{}}(\Qb)^{\lambda({\textbf k{}})} \subset  \\ &
\Mr_{k_1, r^*}(Np, \psi_1;\Qb)\otimes
\Mr_{k_2, r^*}(Np, \psi_2;\Qb)\otimes
\Mr_{k_3, r^*}(Np, \psi_3;\Qb).
\end{align*} 
This map  is then extended to an $\Ar$-linear map
\begin{align}\label{ell}
\ell = 
\ell_{{\textbf f}_1\otimes{\textbf f}_2\otimes{\textbf f}_3, 
\lambda}: \Mc(\Ar)^\lambda\to \Ar;
\end{align}
on the {\sf{} locally free $\Ar$-module of finite rank 
$\Mc(\Ar)^\lambda$}.

This map
produces a sequence of $\Ar$-valued distributions
$\mu^\lambda_r(\chi)\in \Ar$
in such a way that for all suitable weights ${\textbf k{}}\in \Br$ one has
$$
{}
ev_{\textbf k{}}(\mu^\lambda_r(\chi))
=
\scrL( ev_{\textbf k{}}(\pi_\lambda(\Phi_r)(\chi))), \lambda\in \Ar^\times,
\lambda({\textbf k{}})\in \Qb^\times,
{}
$$ 
where ${\textbf k{}}=(k_1{}, k_2{}, k_3{})\in \Br$, $ev_{\textbf k{}}:\Br\to {\mathbb C}_p$
denotes the evaluation map with the property
$$
ev_{\textbf k{}}:\mc(\Ar)\to \mc_{\textbf k{}}({\mathbb C}_p).
$$

More precisely, we consider three
 auxilliary families of modular  forms
{}
\begin{align}\label{EQtifj}&
\tilde  f_{j, k{}_j}(z)=
\\ \nonumber &
\sum_{n=1}^\infty \tilde a_{n,j, k{}_j}e(nz)\in S_{k{}_j}(\Gamma_0(N_jp^{\nu_j}), \psi_j), \ 
(1\le j\le 3, \nu_j\ge 1), 
\end{align}
with the same eigenvalues as those of (\ref{EQfj'}), 
 for all Hecke operators
 $T_q$, with  $q$ prime to $Np$.
In our construction we use as $\tilde f_{j, k{}_j}$ certain
{\sf ``easy transforms''} of  primitive cusp forms 
in (\ref{EQfj}).

In particular, we choose as $\tilde f_j$   certain eigenfunctions
 $\tilde f_{j, k{}_j}=f_{j, k{}_j}^{\,0}$ of 
the  adjoint  Atkin's operator $U_p^*$,
in this case we denote by $f_{j,k{}_j, 0}$    the corresponding eigen\-functions of 
 $U_p$.

The $\Qb$-linear form $\scrL$  produces a ${\mathbb C}_p$-valued admissible measure 
$\tilde \mu^\lambda=\ell( \tilde \Phi^\lambda)$ starting from   the modular $p$-adic admissible 
measure $\tilde \Phi^\lambda$ of stage 3),
where $\ell: {\Mr_T}({\mathbb C}_p)\to{\mathbb C}_p$ denotes a ${\mathbb C}_p$-linear form, 
interpolating $\scrL$. 

\subsection
{Evaluation of $p$-adic integrals}
\subsubsection*{$L$-values and $p$-adic integrals}4)
We show  
that  for any appropriate Dirichlet
character 
$\chi\bmod Np^v$ 
the integral 
$$
{}
\mu^\lambda_{r}(\chi)=\scrL(\pi_\lambda(\Phi_r(\chi)))\in \Ar
$$
evaluated at
$({k{}_1, k{}_2, k{}_3})\in 
\Br= \Br_1\times\Br_2\times\Br_3,  
$
coincides (up to a normalisation) with the special $L$-value  
$$
{}
L^*(f_{1, k{}_1}^{\rho}\otimes f_{2, k{}_2}^{\rho}\otimes f_{3, k{}_3}^{\rho}, 
k{}_2+ k{}_3-2-r, \psi_1\psi_2\chi )
$$
under the above  assumptions on $\chi$ and $r$).
\subsubsection*{A general integral representation 
of Garrett's type} 
The basic idea how a Dirichlet character $\chi$ is incorporated in the 
integral representation \cite{Ga87, BSP96} 
is somewhat similar
to the one used in \cite{B\"o-Sch}, but (surprisingly) more complicated
to carry out.

{\sf Note however
 that the existence of a $\Ar$-valued admissible measure 
$\tilde \mu^\lambda=\ell( \tilde \Phi^\lambda)$ 
established at stage 4), 
does not depend on this technical computation}.

\

In order to control the denominators of the modular forms 
$$
\pi_\lambda(\tilde {\cal E}(-r, \chi))\in{\cal M}^\lambda(\Ar)=:\Phi_r(\chi),
$$used in the construction 
(the admissibility condition)
we use the following result.






 


\section{Criterion of admissibility}

\begin{theo}[Criterion of admissibility ]\label{b3} 
Let $\alpha\in{\cal A}^*$,  $0<|\alpha |_p<1$ and
suppose that 
there exists a positive integer $\varkappa$ such that
the following conditions are satisfied:

1) for all 
$r = 0, 1, {\cdots},  h-1$ 
with $h=[\varkappa  \ord_p\alpha]+1$, and  $v\ge 1$,  
\begin{eqnarray}\label{b6.1}
\Phi _r (a+ (Np^v)) \in
\mc(N p^{\varkappa v}) \qquad\qquad\mbox{\sf{} (the level condition)}
\end{eqnarray}
{}
2) the following 
congruence for the coefficients
holds:
for all 
 $t= 0,1,{\cdots},   h-1$ 
\begin{align}  \label{b6.2}
U^{\varkappa v} \sum^t_{r=0} {t\choose r}
(-a_p)^{t-r} \Phi _r (a+ (Np^v))&\equiv 0 \bmod p^{vt} \\ \nonumber
& \mbox{{\sf{}  (the divisibility condition)}}
\end{align}

Then the linear form given by
$\ds
\tilde \Phi^\alpha(\de_{a+(Np^v)} y^r_p)
:=\pi_\al(\Phi _r(a+(Np^v))
$ 
on local monomials (for all $r=0,1, {\cdots},  h-1$),
is an   $h$-admissible measure:
$
\tilde \Phi^\alpha :{\cal P}^{h  }(Y, \Qb ) \to \mc^\alpha  \subset  
\mc
$
\end{theo}



Proof uses the commutative diagram:
\begin{align*}  
\begin{matrix}
&\mc(Np^{v+1}, \psi; \Ar) &\displaystyle\mathop{\longrightarrow}^ 
{{\pi _{\alpha ,{v}}}}
&\mc{}^\al(Np^{v+1}, \psi; \Ar)& &\cr
&U^{v} \Big\downarrow&&\Big\downarrow\!\!\wr  U^{v}& &\cr
&\mc{}(Np,\psi; \Ar) &
\displaystyle\mathop{\longrightarrow}_ {
{\pi _{\alpha ,{0}}}}
&\mc{}^\al(Np, \psi; \Ar ) &=
\mc{}^\al(Np^{v+1}, \psi; \Ar ).
\end{matrix}
\end{align*}
The existence of the projectors $\pi _{\alpha ,{v}}$
comes from Coleman's Theorem A.4.3 \cite{CoPB}.

\

{\sf On the right: $U$ acts on the locally free $\Ar$-module 
$\mc^\alpha (Np^{v+1}, \Ar)$ }
via the matrix:
\begin{align*}&
\begin{pmatrix} \al &\cdots & \cdots &*\\
0&\al &\cdots &* \\
0 &0 &\ddots &\cdots \\
0&0&\cdots &\al
\end{pmatrix} 
\mbox{ where } 
\al\in \Ar^\times\\ & \Longrightarrow U^v
\mbox{ is an isomorphism  over }\Ar,
\end{align*}
and one controls the denominators of the modular forms
of all levels $v$ by the relation:
\begin{equation}\label{eq:comm}
\pi _{\alpha ,{v}}(h)=U^{-v}\pi _{\alpha ,{0}}(U^vh)=:\pi_\al(h)
\end{equation}
The equality (\ref{eq:comm}) can be used as the definition of 
$\pi_\al$ at any level $Np^v$.

\

The {\sf {} growth condition} 
(see (\ref{equation:adm}))
for $\pi_\al(\Phi_r)$ is deduced 
from the congruences (\ref{b6.2}) between modular forms, using the relation 
(\ref{eq:comm}).
\vfill
{}


\bibliographystyle{amsplain}

\begin{thebibliography}{10xxxx}


\bibitem[Am-V]{Am-V}
{\sc Amice,  Y.} and {\sc Vélu,  J.}, 
{\em Distributions $p$-adiques associ\'ees aux
s\'eries de Hecke}, 
 Journ\'ees Arithm\'etiques de Bordeaux (Conf. Univ.
Bordeaux, 1974), Ast\'e\-risque no. 24/25, Soc. Math. France, Paris
1975, 119 - 131








\bibitem[Boe1]{B\"o1}
{\sc B\"ocherer,  S.}, 
{\em \"Uber die Funk\-tio\-nal\-glei\-chung 
auto\-mor\-pher $L$--Funk\-tio\-nen zur Sie\-gel\-scher Modul\-gruppe}, 
J. reine angew. Math. 362 (1985) 146--168




\bibitem[Boe2]{B\"o2}
{\sc B\"ocherer S.}, 
{\em \"Uber die Fourier--Jacobi Entwickelung
Siegelscher Eisensteinreihen. I.II.}, 
 Math. Z. 183 (1983) 21-46; 189
(1985) 81--100.

\bibitem[BHam]{BHam}  {\sc B\"ocherer, S.},  
{\em Ein Rationalitätssatz für formale Heckereihen zur 
Siegelschen Modulgruppe. }
Abh. Math. Sem. Univ. Hamburg 56, 35--47 (1986)

\bibitem[Boe-Ha]{B\"o-Ha} 
{\sc Boecherer, S., Heim, B.} 
{\em $L$-functions for $GSp_2\times Gl_2$ of mixed weights.}
Math. Z. 235, 11-51(2000)


\bibitem[Boe-Pa6]{Boe-Pa6}
{\sc S. B\"ocherer} and {\sc A.A. Panchishkin},
{\em Admissible $p$-adic measures attached to triple products of elliptic cusp forms},
accepted in Documenta Math. in March 2006 (a special volume  dedicated to John Coates).



\bibitem[BSY]{BSY}
{\sc B\"ocherer, S.},   {\sc Satoh, T.}, and {\sc Yamazaki, T.},
{\em On the pullback of a differential operator and its application to vector valued Eisenstem series, }
Comm. Math. Univ. S. Pauli, 41 (1992), 1-22.




\bibitem[Boe-Schm]{B\"o-Sch}
{\sc B\"ocherer, S.},  and {\sc Schmidt, C.-G.}, 
{\em $p$-adic measures attached to Siegel modular forms}, 
Ann. Inst. Fourier 50, \Numero 5, 1375-1443 (2000).












\bibitem[BoeSP]{BSP96} {\sc B\"ocherer, S.} and {\sc Schulze-Pillot, R.}, 
{\em On the central critical value of the triple product $L$-function. }
David, Sinnou (ed.), Number theory. Séminaire de Théorie des Nombres de 
Paris 1993--94. Cambridge: Cambridge University Press. 
Lond. Math. Soc. Lect. Note Ser. 235, 1-46 (1996). 




\bibitem[Bue]{Bue}  {\sc
Buecker, K.}, {\em 
Congruences between Siegel modular forms on the level of group cohomology.}
Ann. Inst. Fourier (Grenoble) 46 (1996), no. 4, 877--897.



\bibitem[Chand70]{Chand70} 
{\sc Chandrasekharan}, K. Arithmetical functions.
Berlin--Heidelberg--New York,  Springer--Verlag (1970). 


\bibitem[Co]{Co}
{\sc Coates, J.} {\em
On $p$--adic $L$--functions.} Sem. Bourbaki, 40eme annee,
1987-88, n${}^\circ$ 701, Asterisque (1989) 177--178.

\bibitem[Co-PeRi]{Co-PeRi}
{\sc Coates, J.} and {\sc Perrin-Riou, B.}, 
{\em On $p$-adic $L$-functions attached to 
motives over ${\Q}$},
Advanced Studies in Pure Math. 17, 23--54 (1989)



\bibitem[CoPB]{CoPB}
{\sc Coleman, Robert F.}, 
{\em $p$-adic Banach spaces and families of modular forms}, 
Invent. Math. 127, No.3, 417-479 (1997)


\bibitem[CoM]{CoM}
{\sc Coleman, R.} and {\sc Mazur, B.}, 
{\em The eigencurve}, 
Scholl, A. J. (ed.) et al., Galois representations in arithmetic
 algebraic geometry. Proceedings of the symposium,
 Durham, UK, July 9--18, 1996.
 Cambridge: Cambridge University Press. Lond. Math. Soc. Lect. Note Ser. 254,
 1-113 (1998).



\bibitem[CST98]{CST98} {\sc Coleman},  R.,  {\sc Stevens}, G.,  {\sc Teitelbaum}, J., 
        {\it Numerical experiments on families of $p$-adic modular forms}, in:
        Computational perspectives in Number Theory, ed. by 
         D.A. Buell, J.T. Teitelbaum, Amer. Math. Soc., 143-158 (1998).








\bibitem[Colm98]{Colm98}
{\sc Colmez, P.}
{\em Fonctions $L$ $p$-adiques},
 S\'eminaire Bourbaki,
 exposé n${}^\circ$ 851,  novembre 1998.

\bibitem[Colm03]{Colm03}
{\sc Colmez, P.}
La conjecture de Birch et Swinnerton-Dyer $p$-adique. 
Séminaire Bourbaki, 
 exposé n${}^\circ$ 919, juin 2003.

\bibitem[Cour]{Cour} {\sc Courtieu, M.}, 
        {\it Familles d'op\'erateurs sur les formes
modulaires de Siegel et fonctions $L$ $p$-adiques}, \\        
 {\tt http://www-fourier.ujf-grenoble.fr/THESE/ps/ t101.ps.gz},
Th\`ese de Doctorat, Institut Fourier, 2000, 1--122.


\bibitem[CourPa]{CourPa}  {\sc Courtieu,  M.},  {\sc Panchishkin,   A.A.},
          {\it Non-Archimedean $L$-Functions and Arithmetical Siegel Modular Forms},
  Lecture Notes in Mathematics 1471, Springer-Verlag, 2004 (2nd augmented ed.)


\bibitem[De69]{De69}
{\sc Deligne P.}, 
{\em Formes modulaires et représentations $l$-adiques}, 
 Sem. Bourb. 1968/69, exp. no 335. Springer-Verlag, Lect. Notes
 in Math. \Numero 179 (1971) 139 - 172




\bibitem[De79]{De79}
{\sc Deligne P.}, 
{\em Valeurs de fonctions $L$ et 
p\'eriodes d'int\'egrales},
Proc. Symp. Pure Math AMS 33 
(part 2) (1979) 313--342.




\bibitem[Fei86]{Fei86} {\sc Feit, P.}, 
        {\it Poles and Residues of Eisenstein Series for Symplectic
          and Unitary Groups},
        Memoir AMS 61 \Numero 346 (1986), 89 p.


\bibitem[Ga87]{Ga87}
{\sc Garrett,  P.B.}, 
{\em Decomposition of Eisenstein series:
 Rankin triple products}, Ann. of Math. 125 (1987), 209--235


\bibitem[GaHa]{GaHa}
{\sc Garrett,  P.B.} and {\sc Harris, M.}, 
{\em Special values of triple 
 product $L$--functions}, Amer. J. Math 115 (1993) 159--238

\bibitem[Go02]{Go02}  {\sc Gorsse, B.},
Carrés symétriques de formes
modulaires et intégration $p$-adique.
Mémoire de DEA 
de l'Institut Fourier, June 2002


\bibitem[Go-Ro]{Go-Ro}  {\sc Gorsse B.}, {\sc Robert G.},
Computing the Petersson product $\langle f^{\,0},f_0\rangle$.
Prépublication  de l'Institut Fourier, \Numero 654,  (2004).




\bibitem[HaH]{HaH} {\sc Ha Huy, Khoai},
        {\it $p$-adic interpolation and the Mellin-Mazur transform},
        Acta  Math. Viet. 5,  77-99 (1980).


\bibitem[Ka-Za]{Ka-Za}  {\sc Kaneko, M.}, {\sc Zagier, Don,}
 {\it A generalized Jacobi theta function and quasimodular forms.}  
The moduli space of curves (Texel Island, 1994),  165--172, 
Progr. Math., 129, Birkhäuser Boston, Boston, MA, 1995.










\bibitem[Ha93]{Ha93}
{\sc Harris, M.}, and {\sc Kudla, S.},
{\em The central critical value 
 of a triple product $L$-func\-ti\-ons},  
Ann. Math. 133 (1991), 605--672


\bibitem[Hasse]{Hasse} {\sc Hasse, H.},
 {\em Vorlesungen über Zahlentheorie.  }
Zweite neubearbeitete Auflage. 
Die Grundlehren der Mathematischen Wissenschaften, 
Band 59 Springer-Verlag, Berlin-New York 1964 xv+504 pp.


\bibitem[Hi85]{Hi85}
{\sc Hida,  H.}, 
{\em A $p$-adic measure attached to the zeta functions
 associated with two elliptic cusp forms. I}, 
Invent. Math. 79
 (1985) 159--195

\bibitem[Hi86]{Hi86}
{\sc Hida,  H.}, 
{\em Galois representations into $GL_{2}({\Z}_{p}[[X]])$ 
 attached to ordinary cusp forms},
 Invent. Math. 85 (1986) 545--613



\bibitem[Hi90]{Hi90}
{\sc Hida, H.}, 
{\em Le produit de Petersson et de 
Rankin $p$-adique},
S\'eminaire de Th\'eorie des 
Nombres, Paris 1988--1989, 87--102,
Progr. Math., 91, Birkh\"auser 
Boston, Boston, MA, 1990.



\bibitem[Hi93]{Hi93}
{\sc Hida, H.}, 
{\em Elementary theory of $L$-functions and Eisenstein series},
 London Mathematical Society Student Texts. 26 Cambridge:
 Cambridge University Press. ix, 386 p. (1993).






\bibitem[Hi04]{Hi04} H. {\sc Hida},
        {\it $p$-adic automorphic forms on Shimura varieties.}
Springer Monographs in Mathematics.
Springer-Verlag, New York, 2004. xii+390 pp. 


\bibitem[JoH05]{Hu}     {\sc Jory--Hugue, F.}, 
        {\it Unicité des $h$--mesures admissibles $p$-adiques données par des valeurs de
fonctions $L$ sur les caractères}, 
Prépublication de l'Institut Fourier (Grenoble), \Numero 676, 1-33, 2005


\bibitem[Iw]{Iw}
{\sc Iwasawa,  K.}, 
{\em Lectures on $p$--adic $L$--functions},  
Ann. of Math.
 Studies,  \Numero 74. Princeton University Press, 1972


\bibitem[Ibu]{Ibu}
{\sc Ibukiyama, T.} 
{\it On differential operators on automorphic forms 
and invariant pluriharmonic polynomials}, 
Comm. Math. Univ. S. Pauli 48
(1999), 103-118






\bibitem[Ka76]{Ka76}
{\sc Katz, N.M.}, 
{\em $p$-adic interpolation of real analytic Eisenstein
 series}, 
Ann. of Math. 104 (1976) 459--571


\bibitem[Ka78]{Ka78}
{\sc Katz, N.M.}, 
{\em $p$-adic $L$-functions for $CM $--fields}, 
Invent. Math. 48 (1978) 199--297

\bibitem[KiK]{KiK}
{\sc Kitagawa,Koji}, 
{\em On standard $p$-adic $L$-functions of families of
  elliptic cusp forms},
  Mazur, Barry (ed.) et al.: $p$-adic monodromy and the Birch and
  Swin\-ner\-ton-Dyer conjecture. A workshop held August 12-16, 1991 in
  Boston, MA, USA. Providence, R: American Mathematical Society. 
  Contemp. Math. 165, 81-110, 1994


\bibitem[Ku-Le]{Ku-Le}
{\sc Kubota T.} and {\sc Leopoldt H.-W.
}, 
{\em Eine $p$-adische Theorie der
Zetawerte}, 
J. reine angew. Math. 214/215 (1964) 328--339




\bibitem[La]{La}
{\sc Lang S.}, 
{\em Introduction to  Modular Forms}, 
Springer Verlag, 1976








\bibitem[Maa]{Maa71} {\sc Maass, H.},
        {\it Siegel's Modular Forms and Dirichlet Series},
        Springer-Verlag, Lect. Notes in Math. \Numero 216 (1971).



\bibitem[Man73]{Man73}
{\sc Manin,  Yu.I.}, 
{\em Periods of cusp forms and $p$-adic Hecke series},
 Mat. Sbornik 92 (1973) 378--401(in Russian);
Math. USSR, Sb. 21(1973), 371-393 (1975) (English translation). 

\bibitem[Man74]{Man74}
{\sc Manin,  Yu.I.}, 
{\em The values of $p$--adic Hecke series at integer
points of the critical strip}. Mat. Sbornik 93 (1974) 621 - 626
 (in Russian)
 
\bibitem[Man76]{Man76}
{\sc Manin,  Yu.I.}, 
{\em Non--Archimedean integration and $p$-adic
$L$-functions of Ja\-cquet -- Langlands}, 
Uspekhi Mat. Nauk 31 (1976)  5--54
(in Russian);
Russ. Math. Surv. 31, No.1, 5-57 (1976) (English translation). 








\bibitem[Ma-Pa77]{Man-Pa}
{\sc Manin,  Yu.I.} and {\sc Panchishkin,  A.A.
}, 
{\em Convolutions of Hecke
series and their values at integral points}, 
Mat. Sbornik, 104
(1977) 617--651 (in Russian);
Math. USSR, Sb. 33, 539-571 (1977) (English translation).







\bibitem[Ma-Pa05]{Ma-Pa05}
{\sc Manin,  Yu.I.} and {\sc Panchishkin,  A.A.
}, 
{\em Introduction to
Modern Number Theory},
 Encyclopaedia of Mathematical Sciences,  
 vol. 49 (2nd ed.), 
Springer-Verlag, 2005, 514 p.



\bibitem[MTT]{MTT}
{\sc Mazur, B.}, {\sc Tate, J.} and {\sc Teitelbaum, J.}, 
{\em On $p$-adic analogues of the conjectures of Birch and Swinnerton-Dyer}, 
Invent. Math. 84, 1-48 (1986).



\bibitem[Miy]{Miy}
{\sc Miyake, Toshitsune}, 
{\em Modular forms}. Transl. from the Japanese by Yoshitaka Maeda.
Berlin etc.: Springer-Verlag. viii, 335 p. (1989).




\bibitem[Or]{Or}
{\sc Orloff T.}, 
{\em Special values and mixed weight triple products
 (with an Appendix by D.Blasius)}, Invent. Math. 90 (1987) 169--180



\bibitem[Nag1]{Nag1}   {\sc Nagaoka}, S.,
{\it A remark on Serre's example of $p$-adic Eisenstein series},
Math. Z.  235 (2000) 227-250.

\bibitem[Nag2]{Nag2} {\sc Nagaoka}, S.,
{\it Note on ${}\bmod p$ Siegel modular forms},
Math. Z.  235 (2000) 405-420.




\bibitem[Pa94]{Pa94}
{\sc Panchishkin,  A.A.}, 
{\em Admissible Non-Archimedean standard zeta functions
 of Siegel modular forms}, 
 Proceedings of the Joint AMS Summer Conference on Motives,
 Seattle, July 20--August 2 1991, Seattle, Providence, R.I., 1994, vol.2, 
 251 -- 292 


\bibitem[PaViet]{PaViet}
{\sc Panchishkin, A. A.}  
{\em Non-Archimedean Mellin transform and
 $p$-adic $L$-Functions},
Vietnam Journal of Mathematics, 1997, 
N3, 179--202.



\bibitem[PaSE]{PaSE}
{\sc Panchishkin,  A.A.}, 
{\em On the Siegel--Eisenstein measure.
 Israel Journal of Mathematics}, Vol. 120, Part B(2000), 467--509 



\bibitem[PaIAS]{PaIAS}
{\sc Panchishkin,  A.A.}, 
{\em On $p$-adic integration in spaces of modular forms
 and its applications},  
J. Math. Sci., New York 115, No.3, 2357-2377 (2003). 

\bibitem[PaMMJ]{PaMMJ} {\sc Panchishkin,   A.A.},
{\em A new method of constructing $p$-adic $L$-functions associated 
with modular forms}, 
Moscow Mathematical Journal, 2 (2002), Number 
2, 1-16





\bibitem[PaTV]{PaTV} {\sc Panchishkin,   A.A.},
{\em Two variable $p$-adic $L$ functions attached to eigenfamilies of positive slope},
 Invent. Math. v. 154, N3 (2003), pp. 551 - 615






\bibitem[PaB1]{PaB1}
{\sc Panchishkin,  A.A.}, 
{\em Arithmetical differential operators on nearly holomorphic Siegel modular forms},  
Preprint MPI 2002 - 41, 1-52


\bibitem[PaB2]{PaB2}
{\sc Panchishkin,  A.A.}, 
{\em Admissible measures for standard L-functions and nearly holomorphic Siegel modular forms
},  
Preprint MPI 2002 - 42, 1-65

\bibitem[PaJTNB]{PaCS}
 {\sc Panchishkin  A.A.},
{\em Sur une condition suffisante pour l'existence 
des mesures $p$-adiques admissibles}, 
Journal de Théorie des Nombres de Bordeaux, 
v. 15 (2003), pp. 1-24

\bibitem[PaMMJ2]{PaMMJ2}
{\sc Panchishkin,  A.A.}, 
{\em The Maass--Shimura differential operators and 
congruences between 
arithmetical   Siegel modular forms},  
(accepted in Moscow Mathematical Journal  in October 2005 (39 p.).



\bibitem[PSRa]{PSRa}
{\sc Piatetski--Shapiro, I.I.}, and  {\sc Rallis, S.},
{\em Rankin triple
 $L$--functions}, 
 Compositio Math. 64 (1987) 333--399

\bibitem[Puy]{Puy}   {\sc Puydt, J.}, 
        {\it Valeurs spéciales de fonctions $L$ de formes modulaires
 adéliques}, 
 Thèse de Doctorat, Institut Fourier (Grenoble), 19 décembre 2003  

\bibitem[Ran39]{Ran39}
{\sc Rankin, R.A.}, 
{\em Contribution to the theory of Ramanujan's
 function $\tau(n)$ and similar arithmetical functions. I.II},  Proc.
 Camb. Phil. Soc 35 (1939) 351--372



\bibitem[Ran52]{Ran52}
{\sc Rankin, R.A.}, 
{\em The scalar product of modular forms}, 
Proc.
 London math. Soc. 2 (3) (1952) 198-217


\bibitem[Ran90]{Ran90}
{\sc Rankin, R.A.}, 
{\em The adjoint Hecke operator. }
Automorphic functions and their applications, Int. Conf., Khabarovsk/USSR 1988, 163-166
(1990)



\bibitem[Sa]{satoh} {\sc Satoh,T.}, Some remarks on triple L-functions.
Math.Ann.276, 687-698(1987)


\bibitem[Schm]{Schm}
{\sc Schmidt,  C.--G.}, 
{\em The $p$-adic $L$-functions attached to Rankin
 convolutions of modular forms}, J. reine angew. Math. 368 (1986)
201--220




\bibitem[Scho]{Scho}
{\sc Scholl,  A.}, 
{\em Motives for modular forms}, 
Inv. Math. 100 (1990), 419--430


\bibitem[SchEu]{SchEu}
{\sc Scholl, A.J.}, 
{\em An introduction to Kato's Euler systems},
Scholl, A. J. (ed.) et al., Galois representations in arithmetic
 algebraic geometry. Proceedings of the symposium, Durham, UK, July 9-18, 1996.
 Cambridge: Cambridge University Press. Lond. Math. Soc. Lect. Note Ser. 254,
 379-460 (1998).



\bibitem[Se73]{Se73}
{\sc Serre,  J.--P.}, 
{\em Formes modulaires et fonctions z\^eta $p$-adiques}, 
Lect Notes in Math. 350 (1973) 191--268 (Springer Verlag)


\bibitem[SePB]{SePB}
{\sc Serre, J.-P.}, 
{\em  Endomorphisms completement continus des espaces de Banach $p$-adiques}, 
 Publ. Math. Inst. Hautes Etud. Sci., 12, 69-85 (1962). 


\bibitem[Shi71]{Shi71}
{\sc Shimura G.}, 
{\em Introduction to the Arithmetic Theory of
 Automorphic Functions}, Princeton Univ. Press, 1971


\bibitem[Shi75]{Shi75}
{\sc Shimura G.}, 
{\em On the holomorphy of certain Dirichlet series},  
Proc. Lond. Math. Soc. 31 (1975) 79--98


\bibitem[Shi76]{Shi76}
{\sc Shimura G.}, 
{\em The special values of the zeta functions associated with cusp forms.}
  Comm. Pure Appl. Math.  29  (1976), no. 6, 783--804.


\bibitem[Shi77]{Shi77}
{\sc Shimura G.}, 
{\em On the periods of modular forms},  Math. Annalen
229 (1977) 211--221

   
\bibitem[Shi82]{Shi82} {\sc Shimura, G.},
        {\it Confluent Hypergeometric Functions on Tube Domains},
        Math. Ann. 260 (1982), p. 269-302.


\bibitem[Sh83]{Shi83} {\sc Shimura, G.},
        {\it On Eisentsein Series},
        Duke Math. J; 50 (1983), p. 417-476.


\bibitem[ShiAr]{ShiAr}
{\sc Shimura, G.}, 
{\em Arithmeticity in the theory of automorphic forms}, 
Mathematical Surveys and Monographs. 82. Providence, RI:
 American Mathematical Society (AMS) . x, 302 p. (2000)








\bibitem[Til-U]{Ti-U}
{\sc Tilouine, J.} and {\sc Urban, E.
}, 
{\em Several variable $p$-adic
families of Siegel-Hilbert cusp eigenforms 
and their Galois representations}, 
Ann. scient. \'Ec. Norm. Sup. 4${}^{\rm e}$ 
s\'erie, 32 (1999) 499--574.



\bibitem[V]{V}
{\sc Visik, M.M.}, 
{\em Non-Archimedean measures connected with
 Dirichlet series},  
 Math. USSR Sb. 28 (1976), 216-228 (1978). 



\bibitem[MV]{MV}
{\sc Vishik, M.M.} and {\sc Manin, Yu.I.}, 
{\em $p$-adic Hecke series of imaginary quadratic fields},
Math. USSR, Sbornik 24 (1974), 345-371 (1976).



\bibitem[We56]{We56}
{\sc Weil, A.}, 
{\em On a certain type of characters of the 
id\`ele-class group of an algebraic number-field},
Proceedings of the international
symposium on algebraic number theory, 
Tokyo \& Nikko, 1955, pp. 1--7,
Science Council of Japan, Tokyo, 1956.







\bibitem[Wi95]{Wi95}
{\sc Wiles A.}, 
{\em Modular elliptic curves and Fermat's 
Last Theorem},
Ann. Math., II. Ser. 141, No.3 (1995) 443--55.




\end{thebibliography}

\medskip




\end{document}